\documentclass[11pt]{amsart}

\usepackage{amssymb} 
\usepackage[latin1]{inputenc}
\usepackage[matrix,arrow,curve]{xy}
\usepackage[mathscr]{eucal}
\usepackage{epic,eepic}
\usepackage{bbm,url}

\makeindex
\newtheorem{theorem}{Theorem}[section]

\newtheorem{lemma}[theorem]{Lemma}
\newtheorem{prop}[theorem]{Proposition}

\newtheorem{cor}[theorem]{Corollary}

\theoremstyle{definition}

\newtheorem{definition}[theorem]{Definition}
\newtheorem{remark}[theorem]{Remark}

\newcommand{\CC}{\mathcal{C}}

\newcommand{\ord}{\mathrm{ord}}

\newcommand{\Q}{\mathbb{Q}}
\newcommand{\Z}{\mathbb{Z}}
\newcommand{\R}{\mathbb{R}}

\newcommand{\C}{\mathbb{C}}
\newcommand{\W}{\mathcal{W}}

\newcommand{\HH}{\mathcal{H}}

\newcommand{\GL}{\mathrm{GL}}

\newcommand{\Sl}{\mathrm{SL}}
\newcommand{\U}{\mathrm{U}}

\newcommand{\Gl}{\mathrm{GL}}
\newcommand{\Hom}{\mathrm{Hom}}

\newcommand{\N}{\mathbb{N}}

\newcommand{\crit}{\mathrm{crit}}
\newcommand{\coker}{\mathrm{coker\,}}

\newcommand{\SL}{\rm SL}

\renewcommand{\diamond}[1]{\langle #1 \rangle}

\newcommand{\RR}{{\mathcal R}}
\newcommand{\EE}{{\mathcal E}}

\newcommand{\Ac}{{\mathcal A}}
\newcommand{\dual}{\vee}
\newcommand{\hotimes}{\hat \otimes}

\renewcommand{\sp}{\text{Sp\,}}
\newcommand{\anneau}{\mathcal O}

\newcommand{\D}{\mathcal D}

\renewcommand{\AA}{{\mathcal{A}}}

\newcommand{\cont}{\text{cont}}

\newcommand{\cond}{\text{cond}}

\newcommand{\TT}{\mathcal{T}}
\newcommand{\im}{\mathrm{Im\,}}

\newcommand{\Id}{\text{Id}}

\newcommand{\sss}{\text{ss}}

\newcommand{\WW}{\mathcal{W}}

\newcommand{\PP}{{\mathbb{P}}}

\newcommand{\Symb}{{\rm{Symb}}}
\newcommand{\BSymb}{{\rm{BSymb}}}
\renewcommand{\P}{{\mathcal{P}}}

\newcommand{\V}{{\mathcal{V}}}
\renewcommand{\Re}{{\rm{Re}}}

\newcommand{\mat}[1]{\left( \begin{matrix} #1 \end{matrix} \right)}
\newcommand{\DD}{{\mathcal{D}}}

\newcommand{\Mel}{{\rm{Mel}}}
\newcommand{\BCM}{{\rm{BCM}}}
\newcommand{\cusp}{{\rm{cusp}}}
\newcommand{\coim}{{\rm{coker}}}
\newcommand{\Hb}{{\mathbf{H}}}
\newcommand{\res}{{\text{res}}}

\newenvironment{pf}
{\medskip\noindent {\it Proof --- \ }}
{\hfill\nobreak $\Box$ \par\bigbreak}

\begin{document}

\baselineskip 15.8pt

\bibliographystyle{style} 
\title{The $p$-adic $L$-functions of evil Eisenstein series}
\date{}

\begin{abstract} 
We compute the $p$-adic $L$-functions of evil Eisenstein series, showing that they factor as products of two Kubota--Leopoldt $p$-adic $L$-functions times a logarithmic term.
This proves in particular a conjecture of Glenn Stevens.    
\end{abstract}

\author[J.~Bella\"iche]{Jo\"el Bella\"iche}
\address{Jo\"el Bella\"iche\\Brandeis University\\
415 South Street\\Waltham, MA 02454-9110\\U.S.A}
\email{jbellaic@brandeis.edu}
\

\author[S.~Dasgupta]{Samit Dasgupta}
\address{Samit Dasgupta\\ University of California, Santa Cruz\\
1156 High St\\ Department of Mathematics \\ Santa Cruz, CA 95064\\U.S.A}
\email{dasgupta@post.harvard.edu}

\thanks{We thank Glenn Stevens (whose beautiful conjecture was the original motivation for this article) and Robert Pollack
for many useful conversations, encouragement, and access to several of their papers before they were
published. We also thank Henri Darmon who fostered this collaboration.
Jo\"el Bella\"iche was partially supported by NSF grants DMS 0801205 and DMS 1101615, as well as the 2010-2011 
centennial fellowship from the American Mathematical Society.
Samit Dasgupta was partially supported by NSF grants DMS 0900924 and DMS 0952251 (CAREER), as well as a fellowship from the Sloan Foundation.  }
\maketitle

\tableofcontents

\section{Introduction}

The aim of this paper is to compute the $p$-adic $L$-functions of evil Eisenstein 
series (also known as critical Eisenstein series). Before stating our result, 
let us recall how this $p$-adic $L$-function is defined.

Let $f$ be a modular newform of level $\Gamma_1(N)$ and weight $k+2$, with $k\ge 0$ an integer.  Let $p$ be a prime not dividing $N$, and let $\alpha$ and $\beta$ be the roots of the Hecke polynomial $X^2-a_pX+\varepsilon(p) p^{k+1}$, where $\varepsilon$ is the nebentypus of $f$ and $a_p$ is the $T_p$-eigenvalue of $f$.
To attach a $p$-adic $L$-function to $f$,  
one needs  to first choose one of its $p$-refinements $f_\alpha$ or $f_\beta$.
These are forms on \[ \Gamma := \Gamma_1(N) \cap \Gamma_0(p) \]  defined by
  \begin{align*}
f_\alpha(z)&=f(z)-\beta f(pz), \\ 
f_\beta(z)&=f(z)- \alpha f(pz),
\end{align*}
satisfying $U_p f_\alpha = \alpha f_\alpha$, $U_p f_\beta = \beta f_\beta$.
 Choosing one of those two refinements, say $f_\beta$, the $p$-adic $L$-function of $f_\beta$ is traditionally an analytic function $L(f_\beta,\sigma)$ where the variable $\sigma$ runs among continuous characters $\Z_p^\ast \rightarrow \C_p^\ast$.  Already early in the theory, it was observed that the salient $p$-adic object that one attaches to $f_\beta$ is a $p$-adic distribution $\mu_{f_\beta}$ on $\Z_p$, from which we can retrieve the $p$-adic $L$-function by $L_p(f_\beta,\sigma) = \mu_{f_\beta}(\sigma)$.  Here it is understood that $\sigma$ is viewed as a function on $\Z_p$ by extending the character $\sigma$ by $0$ on $p \Z_p$. Note that
the $p$-adic $L$-function determines the restriction of  the  distribution $\mu_{f_\beta}$ to $\Z_p^\ast$,
but ignores the distribution on $p \Z_p$.
Also, it is useful to treat separately the even and odd parts $\mu_{f_\beta}^+$ and $\mu_{f_\beta}^-$ of the distribution $\mu_{f_\beta}$. Each of them determines the values of the $p$-adic $L$-function on half of the characters (the even ones and the odd ones, respectively).

If $\ord_p(\beta) < k+1$, we are in the so-called {\it non-critical case}, and the $p$-adic $L$-function of $f_\beta$ was defined in the 1970s by the work of Mazur and Swinnerton-Dyer, Manin, Visik, and Amice-V\'elu (see e.g.\ \cite{MTT}) by interpolation of the special values of the corresponding Archimedean $L$-function.  To be precise, if $\chi$ is a Dirichlet character of conductor $p^n$ with $n > 0$, and $j$ is an integer in the range $0 \le j \le k$, then the $p$-adic $L$-function satisfies
\begin{equation} \label{noncritinter}
 L_p(f_\beta, \chi z^j) = \frac{p^{n(j+1)} j! }{\beta^n (-2 \pi i)^j G(\chi^{-1})\Omega_f^{\pm}} L(f, \chi^{-1}, j+1).
\end{equation}
(See \cite[Proposition of \S14]{MTT}.)  Here $G(\chi^{-1})$ is the usual Gauss sum, while $\Omega_f^{+}$ and $\Omega_f^{-}$ are
Shimura periods chosen to ensure that the right side of (\ref{noncritinter}) is algebraic; in (\ref{noncritinter}), the plus/minus
sign is determined by $\chi(-1) = \pm 1$.

This definition does not apply in the {\em critical case}, i.e.\ when $\ord_p(\beta) = k+1$.
More recently, Pollack and Stevens provided a definition of $L_p(f_\beta, \sigma)$ when $\ord_p(\beta)=k+1$
 but $f_\beta$ is not in the image the operator $\theta_k$, the so-called {\it non-$\theta$-critical case} (see \cite{stevenspollack2}).  Here $\theta_k$ is the map from overconvergent modular forms of weight $-k$ to overconvergent modular forms of weight $k+2$
that acts as $(q \frac{d}{dq})^{k+1}$ on $q$-expansions. An evil Eisenstein series is $\theta$-critical, so the Pollack--Stevens definition does not apply.  Therefore, we  use the definition of $L_p(f_\beta, \sigma)$ given in \cite{Bcrit} by the first-named author of this paper. This definition extends the Pollack--Stevens construction as follows.

One uses Stevens' notion of {\it overconvergent modular symbols}: these are group homomorphisms from the abelian group $\Delta_0$ of divisors of degree $0$ on the set $\PP^1(\Q)$ to the space of $p$-adic distributions over $\Z_p$, which satisfy a special
$\Gamma$-covariance condition depending on the chosen weight $k$. The space of such  overconvergent modular symbols is denoted  $\Symb_\Gamma(\D_k)$ and is endowed with an action of the traditional Hecke operators, and also an involution $\iota$ that commutes with the Hecke operators. We  denote by $\Symb^\pm_\Gamma(\D_k)[f_\beta]$
the common eigenspace in $\Symb_\Gamma(\D_k)$ for the Hecke operators with the same eigenvalues as $f_\beta$  and for $\iota$ with eigenvalue $\pm 1$.
The main result of \cite{Bcrit} is that, under a mild technical condition on $f$, called {\it decency}, $\Symb_\Gamma^\pm(\D_k)[f_\beta]$ has dimension $1$. We  can  thus, if $f$ is decent, choose generators $\Phi_{f_\beta}^+$ and $\Phi_{f_\beta}^-$ of these spaces, and define the distributions $\mu_{f_\beta}^+$ and 
$\mu_{f_\beta}^-$ as the images of the divisor $\{\infty\}-\{0\}$ under $\Phi_{f_\beta}^+$ and $\Phi_{f_\beta}^-$. We then define the $p$-adic $L$-function by the usual Mellin transform: 
$L_p(f_\beta,\sigma) = \mu_{f_\beta}^\pm(\sigma)$, where the sign $\pm$ is chosen to be $\sigma(-1)$.
Note that since the symbols $\Phi_{f_\beta}^+$ and $\Phi_{f_\beta}^-$ are defined up to multiplication by a non-zero $p$-adic number,
there is the same indeterminacy in the $p$-adic $L$-function; the restriction of $L_p(f_\beta,\sigma)$ to the space of even characters, and to the space of odd characters, are each defined up to multiplication by a non-zero $p$-adic number.

\bigskip

Let us now turn to the case where $f$ is a new Eisenstein series of weight $k+2$ and
level $M$. The complete list of such Eisenstein series is easily given:
there are the {\it normal} ones:
\[ E_{k+2,\psi,\tau}(q)=c_0 + \sum_{n \geq 1} c_n q^n \]
with \[c_n = \sum_{d \mid n} \psi(n/d) \tau(d)d^{k+1}\]
for $n \ge 1$, and \[ c_0 = \begin{cases} 0 & \text{if } Q >1 ,\\ \frac{1}{2}L(\tau, -k-1) =  - \frac{B_{k+2,\tau}}{2(k+2)} & \text{if } Q=1. \end{cases} \]
Here $\psi$ and $\tau$ are primitive Dirichlet characters of conductor $Q$ and $R$ respectively,  such that $QR = M$, $\psi\tau(-1) = (-1)^{k}$,
and if $\psi=\tau=1$, then $k \neq 0$. The latter condition corresponds to the well-known
fact that holomorphic $E_2$ does not exist. Because of this fact, there are also, when $M=\ell$ is prime,
a few {\it exceptional} new Eisenstein series 
$$E_{2,\ell} = \frac{\ell - 1}{24} +  \sum_{n \geq 1}c_n q^n, \text{ with } c_n =\sum_{
\genfrac{}{}{0pt}{}{d|n}{\ell \nmid d}} d^{k+1}.$$
To each new Eisenstein series $f$ as above is attached a sign $\epsilon(f) = \pm 1$, defined as the 
eigenvalue for the $\iota$-involution of the unique (up to scalars) classical modular symbol in $\Symb_{\Gamma_1(M)}(\V_k)[f]$. (See Section~\ref{s:pms} for precise definitions; here $\V_k$ is the dual of the space of polynomials of degree at most $k$.)  One shows (Proposition~\ref{p:boundsymb}) that $\epsilon(f)=\psi(-1)$ if $f=E_{k+2,\psi,\tau}$ and $\epsilon(f)=1$ in the exceptional cases.

If $f$ is a new Eisenstein series as above, and $p \nmid M$, then the two roots $\alpha$ and $\beta$ are $\psi(p)$ and $\tau(p) p^{k+1}$ (resp.\ $1$ and $p^{k+1}$ in the exceptional case).  The form $f_\alpha$ is 
ordinary, whereas the form $f_\beta$ is critical.  The interpolation formula (\ref{noncritinter}) and the 
well-known factorization of the classical $L$-function of $f$ allows one to easily calculate the $p$-adic $L$-function of the ordinary $p$-stabilization $f_\alpha$ in terms of the Kubota-Leopoldt $p$-adic $L$-functions associated to the characters $\psi$ and $\tau$:
\begin{equation} \label{e:ordlp1}
 L_p(f_\alpha, \sigma) = \begin{cases} 0 & \text{if } \sigma(-1) = \epsilon(f) \\
\frac{G(\psi)}{2Q} \sigma^{-1}(Q) L_p(\psi, \sigma z)  L_p(\tau, \sigma z^{-k})& \text{if } \sigma(-1) = -\epsilon(f)
\end{cases}
\end{equation}
in the normal case $f = E_{k+2, \psi, \tau}$, and
\begin{equation} \label{e:ordlp2}
 L_p(f_\alpha, \sigma) = \begin{cases} 0 & \text{if } \sigma(-1) = 1 \\
\frac{1}{2}(1 - \sigma^{-1}(\ell)) \zeta_p(\sigma z) \zeta_p(\sigma) & \text{if } \sigma(-1) = -1
\end{cases}
\end{equation}
in the exceptional case $f = E_{2, \ell}$.
Here the Kubota--Leopoldt $p$-adic $L$-function $L_p(\psi, \sigma)$ is defined by
interpolation of classical values (see equations (\ref{dirichletintera})--(\ref{dirichletinterb}) in the proof of Proposition~\ref{propLpf1} for the precise interpolation formula).  We have chosen the periods $\Omega_f^\pm = -2\pi i$.
A proof of (\ref{e:ordlp1}) is given in Proposition~\ref{propLpf1}, and (\ref{e:ordlp2}) is similar.

Our main result is a similar formula for the critical $p$-stabilization $f_\beta$, when  its $p$-adic $L$-function is defined, i.e.\  when $f$ is decent. For a new Eisenstein series, $f$ is  decent unless  $f$ is a normal Eisenstein series of the form 
$E_{2,\psi,\tau}$ and there exists a prime $\ell$ dividing with the same order $\nu>0$ both the conductors $Q$ and $R$, and such that the restriction of $\psi$ and $\tau$ to $(\Z/\ell^\nu\Z)^\ast$ are equal. In the rest of this article,
we will always assume that our Eisenstein series $f$ is decent.\footnote{The constructions of this paper
do not rely on the assumption that $f$ is decent, i.e.\ even in the indecent case we can construct a partial
modular symbol with the Hecke eigenvalues of $f_\beta$ whose associated $p$-adic $L$-function is given
by Theorem~\ref{t:main}.  However, we do not know in the indecent case that this eigenspace of partial modular symbols has
dimension $1$, nor do we know that the full modular symbol analog
$\Symb^\pm_\Gamma(\D_0)[f_\beta]$ has dimension 1 (which is the reason we cannot define ``the" $p$-adic $L$-function
of $f_\beta$). If we knew both of these one-dimensionality results (for example, if we knew that the eigencurve is smooth
at $f_\beta$), then even in the indecent case $f_\beta$ would have a well-defined $p$-adic $L$-function, and our results
would prove that the formula of Theorem~\ref{t:main} is still valid.}

To state our result, we recall from \cite{BCM} the analytic function $\log_p^{[k]}$, defined as follows. 
Let  $\sigma$ denote a continuous character $\Z_p^\ast \rightarrow \C_p^\ast$.  The function
\[ \frac{d^k \sigma}{dz^k} \cdot \frac{z^k}{\sigma(z)}
\]
is constant on $\Z_p^\ast$; we define $\log_p^{[k]}(\sigma) \in \C_p$ to be this constant.
The analytic function $\log_p^{[k]}$
 vanishes precisely (to order 1) at the characters $\sigma$ of the form $z \mapsto z^j \chi(z)$, where $\chi$ is finite order and $j$ is an integer such that $0 \leq j \leq k-1$.

\begin{theorem} \label{t:main} Let $f$ be a new Eisenstein series of level $M$ and $p$  a prime not dividing $M$. Let $f_\beta$ be the critical slope refinement of $f$. We have
\begin{equation} \label{Lpzero} L_p(f_\beta,\sigma) = 0 \quad \text{ if }\sigma(-1) = - \epsilon(f).\end{equation}
In the normal case $f = E_{k+2,\psi,\tau}$,
we have
\begin{equation} \label{Lpnormalcase}  L_p(f_\beta,\sigma) = \sigma^{-1}(R)\log_p^{[k+1]}(\sigma) L_p(\psi,\sigma z) L_p(\tau,\sigma z^{-k}) \quad \text{ if } \sigma(-1)=\epsilon(f).
\end{equation}
In the exceptional case $f=E_{2,\ell}$, $\ell$ prime, we have
\begin{equation} \label{Lpexpcase} L_p(f_\beta,\sigma) = \log_p^{[1]}(\sigma) (1-\sigma^{-1}(\ell)) \zeta_p(\sigma z) \zeta_p(\sigma) \quad \text{ if } \sigma(-1)=1.
\end{equation}
\end{theorem}
Note that equations (\ref{Lpnormalcase}) and (\ref{Lpexpcase}) are to be interpreted as equalities up to multiplication by a non-zero $p$-adic number, since their left hand sides are
defined only up to multiplication by a non-zero $p$-adic number. However, see Remark~\ref{normalization} below.

\begin{remark} The $-\epsilon(f)$-part of the theorem is relatively easy. The proof is given in Proposition~\ref{p:minuspepsart}.
The $\epsilon(f)$-part is much harder. The proof is given at the end of the paper.
\end{remark}

Let $\chi$ denote a Dirichlet character of $p$-power conductor with $\chi(-1) = \epsilon(f)$.
It is customary to write $p$-adic $L$-functions in terms of a variable $s \in \Z_p$ instead of the variable $\sigma$, writing $L_p(f_\beta,\chi, s)$ for what we call $L_p(f_\beta,\chi\langle z \rangle^{s})$.\footnote{Here we have chosen the convention of \cite{MTT}.  Note that \cite{GS} and other works of Stevens use the convention 
$L_p(f_\beta,\chi, s) = L_p(f_\beta,\chi\langle z \rangle^{s-1})$ exacted by the change of variable $s \mapsto s-1$.}
 Here $\langle z \rangle:= z/\omega(z)$ is the component in $1+p\Z_p$ (or $1+4 \Z_p$ if $p=2$) of $z \in \Z_p^\ast$,  and $\omega$ denotes the Teichm\"uller character.
Similarly, the Kubota--Leopoldt $p$-adic $L$-function is often written as a function of $s \in \Z_p$ via:
\[ L_p(\nu, \chi \langle z \rangle^s) = \begin{cases} L_p(\nu\chi^{-1}\omega, s) & \text{if } \nu\chi^{-1} \text{ is odd} \\
L_p(\nu^{-1}\chi, 1-s) & \text{if } \nu\chi^{-1} \text{ is even} \end{cases}
\]  for any  Dirichlet character  $\chi$  of $p$-power conductor.

With this notation 
equation (\ref{Lpnormalcase}) becomes (for $f=E_{k+2,\psi,\tau}$ and $\chi(-1) = \psi(-1)$)
$$L_p(f_\beta, \chi, s) =  \chi^{-1}(R)\langle R \rangle^{-s} s(s-1)\cdots (s-k)L_p(\psi\chi^{-1},s+1) L_p(\tau^{-1}\chi\omega^{-k},1-s+k).
$$
Equation (\ref{Lpexpcase}) becomes (for $f=E_{2,\ell}$, $\ell$ prime)
\begin{equation} \label{e:exccase}
L_p(f_\beta,s) =  s (1- \langle \ell \rangle^{-s})  \zeta_p(s+1) \zeta_p(1-s).
\end{equation}

In the case $p=3$, $\ell = 11$,  formula (\ref{e:exccase})
was conjectured by Glenn Stevens  based on numerical computations that he carried out with Vincent Pasol using software written by Robert Pollack \cite{stevensEisenstein}.
 
\begin{remark}  In current work,  Ander Steele and Glenn Stevens have taken a different approach towards Theorem~\ref{t:main} (\cite{ss}).
In his thesis, Kalin Kostadinov constructed a $p$-adic family of modular symbols for $\Gamma_0(11)$ and $p=3$ valued in a space of distributions ``with rational poles," and showed that the appropriate specialization of this family yields the formula for the $p$-adic $L$-function  conjectured by Pasol--Stevens \cite{ko}. The work of Steele and Stevens generalizes Kostadinov's result using the Shintani cocycle for $\GL_2(\Q)$  to construct families of modular symbols valued in distributions with poles; they show that this cocycle specializes to the $p$-adic $L$-function of any evil Eisenstein series.
\end{remark}

\begin{remark} Let us indicate some future arithmetic applications of the result of this paper. 

A natural question is whether the main conjecture holds for evil Eisenstein series. A statement of the 
main conjecture in a setting sufficiently general  to contain the case of evil Eisenstein series was first given by Perrin-Riou, and then recently reformulated in terms of $(\phi,\Gamma)$-modules by Pottharst. Thus, by their work, one disposes of an {\it algebraic} $p$-adic $L$-function of an evil Eseinstein series. 
This paper computes a formula for the {\it analytic} $p$-adic $L$-function of an evil Eisenstein series. The main conjecture is the assertion
that the two are equal. This conjecture is proved in a work in preparation by Yurong Zhang, who has computed, using Pottharst's definition and computations of cohomology of $(\phi,\Gamma)$-modules,
a formula for the algebraic $p$-adic $L$-function which matches exactly ours.

Once the main conjecture for evil Eisenstein series is proved, a possible application could be to propagate the main conjecture to classical points on the eigencurve sufficiently close to evil Eisenstein points (or perhaps, to all the classical points on the irreducible components through them), in the spirit of \cite{EPW}. Those points correspond in general to cuspidal non-CM forms with a non-ordinary (but non-critical either) refinement,
and the main conjecture for those forms is not known.

Another application in the same spirit would use the easy part of our theorem, that is, the identical vanishing of our $p$-adic $L$-function
on half of the weight space, to deduce some information about the $\mu$-invariant of ordinary modular forms. This application is an idea
of R. Pollack and a joint work in progress with him.
\end{remark}

\begin{remark} \label{e1} An important open question about the Eisenstein series $f_\beta$
is whether there exist non-classical overconvergent modular forms that are generalized Hecke eigenvectors with the same eigenvalues as $f_\beta$. If we call $e$ the dimension of the space of such forms (classical or not), the question is whether $e=1$.
When $f = E_{k+2,\psi,\tau}$, the following assertions are equivalent.
\begin{itemize}
\item[(a)] We have $e=1$.
\item[(b)] In the category of $p$-adic $G_\Q$-representations, the unique non-split extension of $\psi$
by $\tau(k+1)$ that has good reduction everywhere in the sense of Bloch--Kato is non-split at $p$; that is, the restriction map \[ H^1_f(\Q,\psi^{-1} \tau(k+1))\rightarrow H^1_f(\Q_p, \psi^{-1} \tau(k+1)) \] is injective (hence an isomorphism, since both the source and the target have dimension $1$).
\item[(c)] We have $L_p(\psi^{-1}\tau,z^{k+1}) \neq 0$.
\end{itemize}

The equivalence is proved in \cite{BCsmooth} in the case $\tau=\psi=1$. 
The general case can be proved similarly, using \cite{Bcrit}. It is conjectured that these properties always hold; in fact, (b) is a consequence of Jannsen's conjecture (\cite{jannsen}), cf. \cite[Prediction 5.1]{BKHawaii}. It is widely expected that a proof of such a result would require some progress in transcendence theory (e.g. a suitable generalization of Baker's results  on independence of logarithms, as adapted by Brummer to the $p$-adic setting,  to the case of polylogarithms). 

It is proved in \cite{Bcrit} that if $e \geq 2$, then $L_p(f_\beta,\sigma)$ vanishes at
every {\it interpolation character} $\sigma$, i.e.\  characters of the form $\chi z^j$ with $\chi$ a finite-order character and $j$ an integer $0 \leq j \leq k$. An example is given to show that this does not necessarily hold when $e=1$.

Let us check that this result is compatible with Theorem~\ref{t:main} in the case of a normal Eisenstein series with $\psi=\tau=1$, 
so $k$ is even and $k \geq 2$ (other cases are similar). Observe that the 
factor $\log_p^{[k+1]}(\sigma)$ has a simple zero at every interpolation character and no other zeroes. 
So the only interpolation characters $\sigma$ for which $L_p(f_\beta, \sigma)$ may possibly not vanish are the poles of $L_p(\psi,\sigma z)=\zeta_p(\sigma z)$ and
$L_p(\tau,\sigma z^{-k})=\zeta_p(\sigma z^{-k})$.  The $p$-adic zeta function $\zeta_p(\sigma)$ has (simple) poles precisely at the characters $\sigma(z) = 1$ and $\sigma(z) = z$.
Since we are interested in characters $\sigma$
such that $\sigma(-1)=1$, the only pole of the second factor $\zeta_p(\sigma z^{-k})$ is  $\sigma = z^k$. 
For this $\sigma$, the other term in our factorization formula for $L_p(f_\beta, \sigma)$ is $\zeta_p(\sigma z) = \zeta_p(z^{k+1})$,
which by   the equivalence between (a) and (c) vanishes if and only if $e \neq 1$.
Therefore, we find that $L_p(f_\beta, z^{k})$ is non-zero if and only if $e=1$. 
A similar analysis holds for 
the pole of the first factor $\zeta_p(\sigma z)$ at $\sigma=1$, using the functional equation of the $p$-adic zeta function.  
Our factorization formula therefore conforms with the result proved in \cite{Bcrit} and mentioned in the previous paragraph.
Unfortunately, but not surprisingly, we cannot prove independently that $L_p(f_\beta, z^k) \neq 0$ and therefore conclude the Leopoldt-like conjecture (a) (equivalently, (b) and (c)).

Also, note that in the exceptional case $f=E_{2,\ell}$, one always has $e=1$.  This follows from (\ref{Lpexpcase}), which shows that $L_p(f_\beta,\sigma)$ does not vanish when $\sigma$ is the trivial character.
\end{remark}
\begin{remark} \label{normalization}
With the notation of the preceding remark, assume $e=1$. Then 
it is possible  to reduce the indeterminacy of the modular symbol $\Phi_{f_\beta}^{\epsilon(f)}$, hence of the $p$-adic $L$-function $L_p(f_\beta,\sigma)$. Actually, $e=1$ is equivalent to the fact that the map 
\[ \rho_k : \Symb_{\Gamma}^{\epsilon(f)}(\D_k)[f_\beta] \rightarrow \Symb_{\Gamma}(\V_k)[f_\beta] \] is an isomorphism.
Here, 
$\rho_k$ is induced by the restriction map (i.e.\ $\rho_k$ sends a distribution to the linear form it induces on polynomials of degree at most  $k$). One can then pick an
element $\phi_{f_\beta}$ in $\Symb_{\Gamma}(\V_k)$ and ask that $\rho_k (\Phi_{f_\beta}^{\epsilon(f)})=\phi_{f_\beta}$, 
which reduces the indeterminacy on  $\Phi_{f_\beta}^{\epsilon(f)}$ to whatever indeterminacy we have in our choice of
$\phi_{f_\beta}$. We shall consider two ways to normalize  $\phi_{f_\beta}$.

The first normalization, {\it \`a la Pollack--Stevens}, is the one used in \cite{stevenspollack2} in the cuspidal critical slope case: if $L$ is a finite extension of $\Q_p$ on which $\phi_{f_\beta}$ is defined, and if $\anneau_L$ is its ring of integers,
we normalize $\phi_{f_\beta}$  to be a generator of the free $\anneau_L$-module of rank $1$
$\Symb_{\Gamma}(\V_k(\anneau_L))[f_\beta]$.   This determines $\phi_{f_\beta}$, hence
$\Phi_{f_\beta}^{\epsilon(f)}$ and $L_p(f_\beta,\sigma)$, up to multiplication by an element of $\anneau_L^\ast$ instead of $L^\ast$. In particular, $p$-adic orders of the value of $L_p$ are well-defined. One can then ask for an integral version
of (\ref{Lpnormalcase}), which would hold up to a unit in $\anneau_L$. We shall content ourselves with the case $f=E_{k+2,1,1}$ for an even $k\geq 2$. In this case, assuming $p \neq 2$,
\begin{eqnarray}
\label{lpnorm} L_p(f_\beta,\sigma) = \frac{p^{k-1}}{\zeta_p(k+1) k!} \log_p^{[k+1]}(\sigma) \zeta_p(\sigma z^{-k}) \zeta_p(\sigma z),
\end{eqnarray}
where $L_p(f_\beta,\sigma)$ is normalized \`a  la Pollack-Stevens, and $=$  here  means equality up to a unit in $\anneau_L$.
Note that $\zeta_p(k+1) \neq 0$ since we have assumed $e=1$.
To prove equation~(\ref{lpnorm}),   it suffices to evaluate both sides of (\ref{Lpnormalcase}) at the character $\sigma(z)=z^k$. For the LHS one uses \cite[Example 4.10]{Bcrit}, and for the RHS one uses the well-known residue of $\zeta_p$
at the character $\sigma=z$ (cf. \cite{colmezbourbaki}). It is also possible to prove a formula in the general case using a similar method; we leave this to the interested reader.

The second normalization is the one used by Stevens in the case $f=E_{2,\ell}$ when $\ell=11$ and $p=3$. 
It can be used for the exceptional Eisenstein series $E_{2,\ell}$ without restriction on $\ell$ and $p$ (since $e=1$
in the exceptional case). One simply requires $\phi_{f_\beta}(\{\infty\}-\{0\})=1$, which determines $\phi_{f_\beta}$,
hence $\Phi_{f_\beta}^+$ and $L(f_\beta,\sigma)$, uniquely. One then has 
\begin{equation}   L_p(f_\beta,s) =  \frac{p-1}{p \log_p(\ell)} s (1- \langle \ell \rangle^{-s})  \zeta_p(s+1) \zeta_p(1-s).
\end{equation}
Here $=$ means a genuine equality, and we say that $L_p(f_\beta,s)$ is normalized {\it \`a la Stevens}. The reader may check
that if $p=3$, $\ell=11$, one obtains the precise formula conjectured by Stevens in \cite{stevensEisenstein}.
\end{remark}

We now give some indications about the proof of Theorem~\ref{t:main}. The basic difficulty in computing the $p$-adic $L$-function of a modular form of critical slope is that the interpolation property does not suffice to characterize the $p$-adic $L$-function,
in contrast to the case of non-critical slope (in particular ordinary) modular forms. The basic strategy to overcome this difficulty is the same as the one used in the CM case in \cite{BCM}: there is an injective map 
$$ \Theta_k : \Symb_{\Gamma}(\D_{-2-k}) \longrightarrow \Symb_{\Gamma}(\D_k)$$ induced by the $(k+1)$st derivative on distributions, which serves as
a close analog of the 
aforementioned map $\theta_k$ on overconvergent modular forms. The map $\Theta_k$  commutes with the action of the Hecke operators and the $\iota$ involution, up to a simple twist ($T_\ell$ and $U_\ell$ are multiplied by $\ell^{k+1}$, $\iota$ by $(-1)^{k+1}$). The basic strategy is to find an eigenvector $\Phi_{-2-k}$ in $\Symb_\Gamma(\D_{-2-k})$ whose eigenvalues are such that the
eigenvalues of $\Theta_k \Phi_{-2-k}$ are those of $\Phi_{f_\beta}^\pm$. By uniqueness, we would know 
that $\Theta_k \Phi_{-2-k} = \Phi_{f_\beta}^\pm$ up to a scalar, and from this it is easy to deduce a relation between
the $p$-adic $L$-function of $\Phi_{f_\beta}^\pm$ (the very object we are trying to compute) and the one of 
$\Phi_{-2-k}$. We then hope to get some grasp on the $p$-adic $L$-function of $\Phi_{-2-k}$, which is now an
ordinary modular symbol.  This can perhaps be done by putting it in a family of ordinary modular symbols $\Phi_{k'}$ whose $p$-adic $L$-functions can be computed for positive integer $k'$, or by  other means.

In our case, since $\Phi_{f_\beta}^\pm$ has the eigenvalues (for the $T_\ell$, $U_p$ and diamond operators) of the critical Eisenstein series $f_\beta= E_{k+2,\psi,\tau}^\crit$, the sought-after $\Phi_{-2-k}$ should have the eigenvalues of the ordinary $p$-adic Eisenstein series
$E_{-k,\tau,\psi}^\ord$ (note the negative weight, and the inversion of the order of the characters). It should also
have sign, i.e.\ eigenvalue for $\iota$, equal to $\pm (-1)^{k+1}$. When the sign $\pm$ is $-\epsilon(f_\beta)=-\psi(-1)$,
we therefore want $\Phi_{-2-k}$ to be of sign $\tau(-1)$.
In this case, we can find a suitable $\Phi_{-2-k}$ which is in fact a boundary modular symbol
that we can compute explicitly. It is then easy to compute its $p$-adic $L$-function, which is always $0$. This  allows
the determination of half of the $p$-adic $L$-function of $f_\beta$, and accounts for the easy part of Theorem~\ref{t:main}.

Unfortunately, this method fails for the other, more interesting, part of that $p$-adic $L$-function, 
the one of sign $\pm = \epsilon(f_\beta)=\psi(-1)$. There are simply no modular symbols $\Phi_{-2-k}$ with the right eigenvalues and the right sign $-\tau(-1)$ (at least when $e=1$ in the sense of Remark~\ref{e1}, which is always expected to be the case.)
To solve this problem, we employ the notion of {\it partial modular symbols} due to the second-named author (see \cite{D}, \cite{DD}). A partial modular symbol is a modular symbol defined only on the divisors of degree $0$ on a non-empty subset $C$
of the set of cusps $\PP^1(\Q)$, instead of the whole set of cusps $\PP^1(\Q)$ as for a usual {\it full} modular symbol.

From the point of view of partial modular symbols for a subset $C$ of cusps, an Eisenstein series that vanishes at all the cusps in $C$
looks like a cusp form.  It is therefore possible to attach to such an Eisenstein series a partial modular symbol by integration, exactly as one does in the Manin--Shokurov construction of modular symbols attached to  cusp forms. For example, if the Eisenstein series $E_{k+2,\tau,\psi}$ vanishes at both the cusps $0$ and 
$\infty$ (this is the case when both $\tau$ and $\psi$ are non-trivial), it is possible
to attach by integration (and Stevens' lifting) an ordinary  partial modular symbol $\Phi_{k}$ in $\Symb_{\Gamma,C}(\D_k)$ for a suitable set of cusps $C$ containing $0$ and $\infty$ with the same eigenvalues as $E_{k+2,\tau,\psi}^\ord$ and {\it the right sign $-\tau(-1)$}.

Moreover, the very construction by integration of $\Phi_k$ allows us to compute its $p$-adic $L$-function
as a product of two Dirichlet $p$-adic $L$-functions. By interpolation, using computations made in \cite{DD}, we can define
a partial modular symbol $\Phi_{-2-k}$ in $\Symb_{\Gamma,C}(\D_{-2-k})$  with the right eigenvalues and sign.

To finish the proof we need to deduce that $\Theta_k \Phi_{-2-k}$ is the same, up to a non-zero scalar, as the restriction to the set of cusps $C$ of the full modular symbol $\Phi_{f_\beta}^{\psi(-1)}$. For this we need two things: to know that the eigenspace we are considering in the
space of {\it partial} modular symbols has dimension $1$, and that the restriction of  $\Phi_{f_\beta}^{\psi(-1)}$ to a partial modular symbol over $C$ is still non-zero. The first of these facts is proven with a method that is similar to the one used to \cite{Bcrit} to prove the same result for full modular symbols: it involves constructing the eigencurve for partial modular symbols and comparing it with various other eigencurves. The second of these facts is perhaps the most technically difficult point of the paper: it involves a very careful
study of boundary overconvergent modular symbols. 

This description of our method applies for an Eisenstein series that vanishes at $0$ and $\infty$, but not all Eisenstein series are of this type. However, if $f$ is an Eisenstein series, it is always possible to chose two auxiliary primes $\ell_1$ and $\ell_2$
and a linear combination $g$ of $f(z),\  f(\ell_1z),\  f(\ell_2z)$ and $f(\ell_1 \ell_2 z)$ that vanishes at both $0$ and $\infty$. By applying the method described above to $g$, with some complication due to the fact that $g$ is no longer a newform,
we eventually get a formula  for the $p$-adic $L$-function of $f_\beta$ which is (\ref{Lpnormalcase})  up to parasitic factors involving $\ell_1$ and $\ell_2$. An important point in the proof is the fact that the eigencurve is still smooth at the old point $g$, which can be proved assuming that $\ell_1$ and $\ell_2$ are outside a set of {\it bad} primes depending on $f$. When $f$ is a normal Eisenstein series, the set of bad prime is finite, and  by letting the auxiliary primes $\ell_1$ and $\ell_2$ vary, we can show that the parasitic factors cancel, yielding our desired result in the normal case. This methods breaks down
in the exceptional case $f=E_{2,\ell}$ because for all choices of an auxiliary prime $\ell_1$, the eigencurve of tame level $\ell \ell_1$ is non-smooth (even non-irreducible) at the old point corresponding to $f$. We therefore need an alternative method for the final steps of our proof, which we provide in the last section, based on a certain numerical coincidence that arises 
only in the exceptional case.

\section{Partial modular symbols} \label{s:pms}

\subsection{Hecke operators}

\label{hecke} 
Let  $\Gamma$ and $\Gamma'$ be congruence subgroups of $\Gl_2(\Q)$, $S$ a submonoid of $\Gl_2(\Q)$ containing $\Gamma$ and $\Gamma'$,  and $W$ a right $S$-module.  For each double coset $\Gamma s \Gamma'$ with $s \in S$, we have a morphism 
\[  [\Gamma s \Gamma']_W: W^\Gamma \rightarrow W^{\Gamma'} \] called a Hecke operator, defined by \[ w_{|[\Gamma s \Gamma']} = \sum_i w_{|s_i}.\] Here the $s_i \in S$ are defined by \[ \Gamma s \Gamma' = \coprod_i \Gamma s_i  \text{ (finite decomposition)}.\] We often drop  $W$ from the notation of the Hecke operators.

In this paper,  $W$ will always be a vector space over a field $L$ of characteristic $0$, and the action of $S$ will be $L$-linear.
It follows that the Hecke operators are also $L$-linear.

We recall standard names for some Hecke operators that we will use throughout the paper. 
When $\Gamma=\Gamma'$ is $\Gamma_1(N)$,  $\Gamma_0(N)$, or any group in between, we denote by
$T_\ell$ (resp. $U_\ell$) the Hecke operator 
$[\Gamma \mat{1 & 0 \\ 0 & \ell} \Gamma]$ for $\ell \nmid N$ (resp. for $\ell \mid N$).
Note that $$\Gamma \mat{1 & 0 \\ 0 & \ell } \Gamma = \coprod_{a=0}^{\ell-1} \Gamma \mat{1 & a \\ 0 & \ell } \coprod  \Gamma \mat{\ell & 0 \\ 0 & 1}$$ if $\ell \nmid N$.  The same holds without the last coset for $\ell \mid N$.
For $a \in (\Z/N\Z)^\ast$, we denote by $\langle a \rangle$ the Hecke operator $\Gamma s_a \Gamma$, where $s_a$ is any matrix in $\Gamma_0(N)$ whose upper-left entry is congruent to $a \pmod{N}$. All these operators commute.

Now let $M$, $N$, and $t$ be positive integers such that $Mt \mid N$.  Let $\Gamma = \Gamma_1(M)$ or $\Gamma_0(M)$, and 
let $\Gamma' = \Gamma \cap \Gamma_0(N)$.
Let $\alpha_t =  \mat{1 & 0 \\ 0 & t^{-1}}$
and denote by $V_t$ the operator \[ t^{-1} [\Gamma \alpha_t \Gamma']: W^{\Gamma} \rightarrow W^{\Gamma'}. \] Observe that 
$\Gamma\alpha_t \Gamma' = \Gamma \alpha_t,$
so $w_{|V_t}= t^{-1} w_{|\alpha_t}.$  It follows that $V_{tt'} = V_t V_{t'}$ whenever 
$Mtt' \mid N$,
so it is enough to consider $V_\ell$ for $\ell$ a prime factor of $N/M$. It is clear that $V_\ell$ commutes with any operator $T_{\ell'}$ or $U_{\ell'}$ for $\ell' \neq \ell$, as well as with the diamond operators.

A simple computation shows that if we write  $i_{\Gamma,\Gamma'}$ for the inclusion
 $W^{\Gamma} \rightarrow W^{\Gamma'}$, then
\begin{equation} \label{UlVl}  w_{|V_\ell  U_\ell}  = i_{\Gamma,\Gamma'}(w), \end{equation}
and if $\ell \mid N$ but $\ell \nmid M$, then
\begin{equation} \label{UlVlTl} i_{\Gamma,\Gamma'}(w_{|T_\ell}) =  i_{\Gamma,\Gamma'}(w)_{|U_\ell} +  \ell w_{|V_\ell T(\ell,\ell)},
\end{equation}  
where $T(\ell,\ell)$ is the operator $\Gamma' \mat{\ell & 0 \\ 0 & \ell} \Gamma'$. 

Finally, we let $\alpha_\infty:=\left( \begin{matrix} 1 & 0 \\ 0 & -1 \end{matrix} \right) \in S$ and define an involution
$\iota\colon W^\Gamma \rightarrow W^\Gamma$ by
 $ \iota=\left[ \Gamma \alpha_\infty \Gamma \right]. $
Note that $\iota$ is simply given by $w \mapsto w_{|\alpha_\infty}.$

Let $W_{\text{univ}}$ be the $\Q$-vector space of maps  $S \rightarrow \Q$, endowed with its obvious right $S$-action by left-translations. 
We write $\HH(S,\Gamma)$ for the $\Q$-algebra generated by all Hecke operators $[\Gamma s \Gamma]$ acting on $W_{\text{univ}}$.
It is easy to see that $\HH(S,\Gamma)$ acts on every right $S$-module $W$ (by letting $[\Gamma s \Gamma]_{W_\text{univ}}$ act by $[\Gamma s \Gamma]_W$), making $W$ a right $\HH(S,\Gamma)$-module.

The basic example of these actions is when $W$ is the space of holomorphic functions $f\colon \Hb \rightarrow \C$, where $\Hb$ is the Poincar\'e upper half-plane,
and $S=\Gl_2^+(\Q)$ (matrices with positive determinant) acts by
\begin{equation} \label{actionheckemodular} f_{|\gamma}(z) = (\det\gamma)^{k+1}  (cz+d)^{-k-2} f\left( \frac{az+b}{cz+d}\right) \text{ for }\gamma=\mat{a & b \\ c &d},\end{equation}  for a fixed integer $k$. 
Then $W^\Gamma$ contains the space of modular forms of level $\Gamma$ and weight $k+2$
 as a subspace stabilized by the action of $\HH(\Gl_2^+(\Q),\Gamma)$. The reader can check
 that with these conventions, if $f(z)=\sum_{n=0}^\infty a_n q^n$ with $q = e^{2 i \pi z}$, then 
 \[ f_{|V_\ell}(z) = f(\ell z) = \sum_{n=0}^\infty a_n q^{\ell n} \] and \[ f_{|U_\ell}(z)=\sum_{n=0}^\infty a_{n \ell} q^{nz}. \]

\subsection{Definition of partial modular symbols} \label{s:dpms}

We let $\Gl_2(\Q)$ act on the left on $\PP^1(\Q)$ by $\gamma \cdot x=\frac{ ax+b}{cx+d}$.
Let $\Gamma$ be a congruence subgroup of $\SL_2(\Z)$. Let $C$ be a non-empty $\Gamma$-invariant subset of $\PP^1(\Q)$. We  denote by $\Delta_C$ the abelian group of divisors on $C$, i.e.\  \[ \Delta_C =\left\{ \sum_{c \in C} n_c \{c\}: \ n_c \in \Z,\ \ n_c=0 \text{ for almost all }c\right\}, \] and by $\Delta_C^0$ the subgroup of divisors of degree $0$ (i.e., such that $\sum_{c \in C} n_c=0$). The group $\Gamma$ acts naturally on the left on $\Delta_C$ and $\Delta_C^0$.

If $V$ is any abelian group endowed with a right-action of $\Gamma$, we provide
$\Hom(\Delta_C^0,V)$ with a right $\Gamma$-action by setting, for 
$\phi \in \Hom(\Delta_C^0,V)$, \[ \phi_{|\gamma}(D) = \phi(\gamma \cdot D)_{|\gamma}.\]

\begin{definition} A {\it partial modular symbol} on $C$ for $\Gamma$ 
with values in $V$ is a $\Gamma$-invariant element $\phi \in  \Hom(\Delta_C^0,V)$. We write \[ \Symb_{\Gamma,C}(V) := \Hom(\Delta_C^0,V)^\Gamma \] for
the abelian group of partial modular symbols. 

When $C = \PP^1(\Q)$, we drop $C$ from the notation and  
call $\Symb_\Gamma(V)$ the space of {\it (full) modular symbols} for $\Gamma$ 
with values in $V$.
\end{definition}

Assume that we are given a submonoid $S$ of $\Gl_2(\Q)$ containing $\Gamma$
that acts on the right on $V$ in a way extending the action of $\Gamma$,
and that preserves the set $C \subset \PP^1(\Q)$. Then $\Hom(\Delta_C^0,V)$
has a natural right action of $S$ extending the one of $\Gamma$, and thus
$\Symb_{\Gamma,C}(V)$ is endowed with a right action of the Hecke algebra $\HH(S,\Gamma)$.
It is clear that $\Symb_{\Gamma,C}(V)$ is a functor from the category of right $\Gamma$-modules $V$ (resp.\ right $S$-modules $V$) to the category of abelian groups (resp.\ right $\HH(S, \Gamma)$-modules).  
 
\subsection{Comparison between partial and full modular symbols}

The exact sequence of abelian groups $$0 \rightarrow \Delta_C^0 \rightarrow \Delta_{C} \rightarrow \Z \rightarrow 0$$ is split, and hence we have an exact sequence
of $S$-modules:
\begin{equation} \label{exact1}
 0 \rightarrow V \rightarrow \Hom(\Delta_C,V) \rightarrow \Hom(\Delta_C^0,V) \rightarrow 0.
\end{equation}
The long exact sequence of group cohomology for $\Gamma$ attached to this short exact sequence gives an exact sequence of $\HH(S,\Gamma)$-modules:
\begin{equation} \label{exact2}
\xymatrix{ 0 \ar[r] & V^\Gamma \ar[r] & \Hom(\Delta_C,V)^\Gamma \ar[r]^{\ b_C} &
 \Symb_{\Gamma,C}(C)\ar[r]^{\ h_C} & H^1(\Gamma,V)}.
\end{equation}
\begin{definition} We write $\BSymb_{\Gamma,C}(V)$ for the image of the map $b_C$ in the above exact sequence, and call this submodule of $\Symb_{\Gamma,C}(V)$ the module of {\it boundary partial modular symbols}.
\end{definition} 
By (\ref{exact2}), we have an isomorphism of  $\HH(S,\Gamma)$-modules
\begin{equation} \label{bsymbvg}\BSymb_{\Gamma,C}(V) \simeq  \Hom(\Delta_C,V)^\Gamma /V^\Gamma \end{equation}
and an exact sequence of $\HH(S,\Gamma)$-modules
\begin{equation} \label{exact3}
\xymatrix{ 0 \ar[r] & \BSymb_{\Gamma,C}(V) \ar[r] & 
 \Symb_{\Gamma,C}(V) \ar[r]^{h_C} & H^1(\Gamma,V).}
\end{equation}

When $C=\PP^1(\Q)$, this exact sequence is well-known and we have a commutative diagram:
$$\xymatrix{ 0 \ar[r] & \BSymb_{\Gamma}(V) \ar[r] \ar[d] &
 \Symb_{\Gamma}(V) \ar[r] \ar[d]^{\res_C} &  H^1(\Gamma,V) \ar[d]^{\simeq}  
 \\
0 \ar[r] & \BSymb_{\Gamma,C}(V) \ar[r] &
 \Symb_{\Gamma,C}(V) \ar[r] &  H^1(\Gamma,V) }$$
 
An easy diagram-chase yields:
 \begin{lemma} \label{kerrestriction} The kernel of the restriction map 
 \[\res_C\colon \Symb_{\Gamma}(V) \rightarrow \Symb_{\Gamma,C}(V) \] is contained in $\BSymb_{\Gamma}(V).$
 \end{lemma}

\subsection{Cohomological interpretation}

In this subsection we mimic the arguments of \cite[Prop 4.2]{ashstevens} to give a
cohomological interpretation of partial modular symbols.

Let $\Hb$ be the Poincar\'e upper half-plane, and $\Hb_C = \Hb \sqcup C$ with its usual topology. Let $V$ be the constant sheaf on $\Hb_C$. The long exact sequence of relative cohomology  of $(\Hb_C,C)$ gives
$$0 \rightarrow H^0(\Hb_C,C,V) \rightarrow H^0(\Hb_C,V) \rightarrow H^0(C,V) \rightarrow H^1(\Hb_C,C,V) \rightarrow H^1(\Hb_C,V).$$
The first term $H^0(\Hb_C,C,V)$ vanishes because $\Hb_C$ is connected and $C$ is not empty. The last term $ H^1(\Hb_C,V)$ also vanishes because $\Hb_C$ is contractible. We therefore obtain an exact sequence  of $\Gamma$-modules
\begin{equation} 
0 \rightarrow H^0(\Hb_C,V) \rightarrow H^0(C,V) \rightarrow H^1(\Hb_C,C,V) \rightarrow 0. \label{e:hcohom}
 \end{equation} 
We claim that (\ref{e:hcohom})
 is isomorphic as an exact sequence of $S$-modules to the sequence (\ref{exact1}). Indeed there is a natural isomorphism
$H^0(\Hb_C,V) \rightarrow V$ since $\Hb_C$ is connected, and a natural isomorphism $H^0(C,V) \rightarrow \Hom(\Delta_C,V)$ since $C$ is discrete. Thus there is an isomorphism $H^1(\Hb_C,C,V) \rightarrow \Hom(\Delta_C^0,V)$. By taking the $\Gamma$-invariants, there is an isomorphism  
$H^1(\Hb_C,C,V)^\Gamma \rightarrow \Symb_{\Gamma,C}(V)$. The Hochschild-Serre spectral sequence then gives:
\begin{prop} \label{propisomcohom} There is a natural $\HH(S,\Gamma)$-equivariant isomorphism, functorial in $V$:
$$\Symb_{\Gamma,C}(V) \cong H^1(\Gamma \backslash \Hb_C,
\Gamma \backslash C,\tilde V),$$
where $\tilde V$ is the sheaf on $\Gamma \backslash \Hb_C$ corresponding to $V$.
\end{prop}

\subsection{Classical partial modular symbols}

\label{classicalpartial}

Let $L$ be a field of characteristic $0$ and $k \geq 0$ an integer. 
Recall that the algebra $\HH(\GL_2^+(\Q), \Gamma)$ acts naturally on the space of modular forms
$M_{k+2}(\Gamma,L)$, cusp forms $S_{k+2}(\Gamma,L)$, and Eisenstein series
$\EE_{k+2}(\Gamma,L)$.
 Let $\P_k(L)$ be the 
$L$-space of polynomials in one variable $z$ of degree at most $k$. 
Let $S=\Gl_2(\Q)$ act on the left on $\P_k(L)$ by
\begin{eqnarray} \label{actpk} (\gamma \cdot P)(z)
=	(a-cz)^k P\left(\frac{dz-b}{a-cz}\right).
\end{eqnarray}
Define $\V_k(L) = \P_k(L)^\dual = \Hom_L(\P_k,L)$ with a right action
of $S$ given by
$$f_{|\gamma}(P) = f(\gamma \cdot P)  \quad \text{for all } P \in \P_k(L).$$
The space $\Symb_\Gamma(\V_k(L))$ is called the space of {\it classical modular symbols} of weight $k$ over $L$. 
We call $\Symb_{\Gamma,C}(\V_k(L))$ the space of {\it classical partial modular symbols}.

We write $h$ for the natural map
 $$h:\Symb_\Gamma(\V_k(L)) = H^1_c (\Gamma,\V_k(L)) \rightarrow H^1(\Gamma,\V_k(L)).$$ 
The following proposition is a version of the classical Eichler-Shimura isomorphism:
\begin{prop} \label{p:bsymb}
Assume that the congruence subgroup $\Gamma$ satisfies $\Gamma_1(N) \subset \Gamma \subset \Gamma_0(N)$ for some integer $N$. Then, after possibly replacing $L$ by a finite extension, the following holds:
\begin{itemize}
\item[(i)]
We have an exact sequence of $\HH(\Gl_2^+(\Q),\Gamma)$-modules 
$$ 0 \longrightarrow \BSymb_{\Gamma}(\V_k(L)) \longrightarrow  \Symb_{\Gamma}(\V_k(L)) 
 \stackrel{h}{\longrightarrow} H^1(\Gamma,\V_k(L)) \longrightarrow \EE_{k+2}(\Gamma,L) \longrightarrow 0.$$
\item[(ii)]
There exists an isomorphism of $\HH(\Gl_2^+(\Q),\Gamma)$-modules
$$ \im h \simeq S_{k+2}(\Gamma,L)^2.$$
\item[(iii)] 
There exists an isomorphism 
$$ \BSymb_{\Gamma}(\V_k(L)) \simeq \EE_{k+2}(\Gamma,L) $$
that is compatible with the Hecke operators $T_\ell$ for $\ell \nmid N$, $U_\ell$ for $\ell \mid N$, and $\langle a \rangle$
for $a \in (\Z/N\Z)^\ast$.
\end{itemize}
\end{prop}
\begin{pf} 
The point (ii) is the most classical. A convenient reference is Theorem 1 of \cite[\S6.3]{hidabb}. Let us recall
the salient point of the proof.  One defines an $\R$-linear map compatible with the action of
$\HH(\Gl_2^+(\Q),\Gamma)$, 
\begin{align*}
 \phi : S_{k+2}(\Gamma,\C) &\rightarrow \Symb_{\Gamma}(\V_k(\R))  \\
 f & \mapsto \phi_f,
 \end{align*} by the formula \[ \phi_f(\{x\}-\{y\})(P)=\Re \int_y^x f(z) P(z) dz \] for any two cusps $x$ and $y$. 
One proves that $h \circ \phi$ is injective by proving that it transforms the real or imaginary part  of the Peterson inner product 
(in the cases $k$ is odd or $k$ is even, respectively) into the Poincar\'e duality product. So $h \circ \phi$ defines an injective map
$ S_{k+2}(\Gamma,\C) \rightarrow \im h$ which is proven to be an isomorphism by equality of dimension of the 
target and the source. Complexifying, we get the desired isomorphism over $\C$. To deduce it over
some finite extension of $\Q$, we need to choose a period for each new form.

To prove (i), we construct a map $\psi :\EE_{k+2}(\Gamma,\C)  \rightarrow H^1(\Gamma,\V_k(\C))$ by sending $f$ to
the class of the cocycle \[ \gamma \mapsto \left(P \mapsto \int_{x}^{\gamma x} f(z)P(z) dz\right), \] where $x$ is a fixed point 
in $\Hb$. It is clear that this map does not depend on $x$, is injective, and is compatible with the action of $\HH(\Gl_2^+(\Q),\Gamma)$.
The image of $\psi$ has trivial intersection with $\im h$ since Eisenstein series and cupsidal forms have different systems of Hecke eigenvalues, and actually $H^1(\Gamma,\V_k(\C)) = \im h \oplus \im \psi$ by dimension-counting. The result follows.

To show (iii), note that $\BSymb_\Gamma(\V_k(L)) = \ker h$ is in perfect duality
(with respect to the  Poincar\'e pairing) with $\coim \ h = \EE_{k+2}(\Gamma,L)$. For this duality the $T_\ell$ are self-adjoint, but not the $U_\ell$ and $\langle a \rangle$. However, there exists a slight modification
 of Poincare's pairing (\cite[exercise 5.5.1]{diamond}) that is still perfect and for which the $U_\ell$ and 
 $\langle a \rangle$ as well as the $T_\ell$ are self-adjoint. This gives an isomorphism $\BSymb_\Gamma(\V_k(L)) \simeq
\EE_{k+2}(\Gamma,L)^\dual$. Since it is not hard to prove (cf. e.g. \cite{Bcourse}) that  $\EE_{k+2}(\Gamma,L)$ is self-dual as a module over the Hecke operators $T_\ell$, $U_\ell$, and $\langle a \rangle$, the result follows. 
\end{pf}

\begin{remark} We warn the reader that if $M, N$, and $t$ satisfy $Mt \mid N$, then the isomorphisms of Proposition~\ref{p:bsymb}(iii) for $\Gamma = \Gamma_1(M), \Gamma_1(N)$ may not in general be chosen to be compatible with 
the Hecke operators $V_t$. 
Actually, in \S\ref{ss:bs}, we shall provide an independent
proof of the special case of (iii) that we need in this paper. Namely, for $f$ a Eisenstein series for $\Gamma_1(M)$, and $\Gamma=\Gamma_1(N)$ for $N$ a multiple of $M$ satisfying certain conditions (cf. equation (\ref{a:sf})), we shall construct an explicit isomorphism 
\[ \BSymb_\Gamma(\V_k(L))[f] \simeq \EE_{k+2}(\Gamma,L)[f] \] 
and we shall also compute the actions of the operators $V_t$ on both sides. This computation will play an essential role in section
5, and  will also make clear that in some cases no isomorphism can respect these operators.
\end{remark}

\begin{cor} \label{corparsymb} Retain the assumptions of  Proposition~\ref{p:bsymb}. Let $C$ be non-empty set of cusps that is stable under 
$\Gamma$. As a module over the algebra generated by the Hecke operators $T_\ell$ (or $U_\ell$) for all primes $\ell$ such that $\mat{1 & 0 \\ 0 & \ell}$ stabilizes $C$, and $\langle a \rangle$ for all $a \in (\Z/N\Z)^\ast$ such that $s_a$
stabilizes $C$, the semi-simplification of $\Symb_{\Gamma, C}(\V_k(L))$ is isomorphic to a
submodule of two copies of $M_{k+2}(\Gamma,L)$. \end{cor}
\begin{pf} 
We have the commutative diagram
$$\xymatrix{ 0 \ar[r] & \BSymb_\Gamma(\V_k(L)) \ar[r] \ar[d]^{a} &
 \Symb_{\Gamma}(\V_k(L)) \ar[r]^h \ar[d] &  H^1(\Gamma,\V_k(L)) \ar[r] \ar[d]^{=}  
&  \EE_{k+2}(\Gamma,L) \ar[r] \ar[d]^{b} & 0 \\
0 \ar[r] & \ker h_C \ar[r] &
 \Symb_{\Gamma,C}(\V_k(L)) \ar[r]^{h_C} &  H^1(\Gamma,\V_k(L)) \ar[r] & \coker h_C \ar[r] & 0,  }$$
 from which we get two exact sequences (compatible with the Hecke operators listed in the statement of the corollary):
 $$0 \longrightarrow \ker h_C \longrightarrow \Symb_{\Gamma,C}(\V_k(L)) \longrightarrow \im h_C \longrightarrow 0$$
 and
 $$0 \longrightarrow \im h \longrightarrow \im h_C \longrightarrow \ker b \longrightarrow 0.$$
 After semisimplification, we thus obtain an injective morphism 
 $$\Symb_{\Gamma,C}(\V_k(L))^\sss \subset (\im h)^\sss \oplus (\ker h_C)^\sss \oplus (\ker b)^\sss.$$
 By Proposition~\ref{p:bsymb}, $(\im h)^\sss \subset (S_{k+2}(\Gamma,L)^\sss)^2$.  The leftmost vertical arrow $a$ is identified with 
 \[ \BSymb_{\Gamma}(\V_k(L)) \rightarrow \BSymb_{\Gamma, C}(\V_k(L)), \] 
which is clearly surjective.  Hence $\ker h_C$ is a quotient of $\EE_{k+2}(\Gamma)$ according to the proposition, yielding an injection $(\ker h_C)^\sss \subset  \EE_{k+2}(\Gamma,L)^\sss$.   Clearly $\ker b \subset \EE_{k+2}(\Gamma,L)$.  The result follows.
 \end{pf}

The proof of Corollary~\ref{corparsymb} explains our interest in partial modular symbols.  By restricting to a subset $C$ of the cusps we remove some boundary symbols, but more importantly we gain some ``honest" modular symbols associated to Eisenstein series (namely those Eisenstein series that are $C$-cuspidal); these symbols are defined by periods on the upper half plane and are hence related to values of classical $L$-functions.  This is explored in detail in \S\ref{s:eismodsymb}.

\subsection{The boundary symbol associated to an Eisenstein series} \label{s:bound}

Consider the new Eisenstein series $E_{k+2,\psi,\tau}$ and $E_{2,\ell}$ defined in the introduction. If $t$ is a positive integer, we write
 \[ E_{k+2,\psi,\tau,t}=(E_{k+2,\psi,\tau})_{|V_t}\] 
 (i.e.\ $E_{k+2,\psi,\tau,t}(z)=E_{k+2,\psi,\tau}(tz)$) and $E_{2,\ell,t}=(E_{2,\ell})_{|V_t}$. We recall the following classification of Eisenstein series, which can be found in \cite{miyake} and in the precise form given below in \cite{stein}.
\begin{prop}\label{p:miyakestein} 
 Let $L$ be a field of charcteristic $0$ containing the $\varphi(N)$th roots of unity.
 For $k>0$, 
the series $E_{k+2,\psi,\tau,t}$ for $\psi$ a Dirichlet character of conductor $Q$ and
$\tau$ a Dirichlet character of conductor $R$ such that $\tau\psi(-1)=(-1)^k$, and $t$ a positive integer such that $QRt|N$, form a basis of $\EE_{k+2}(\Gamma_1(N),L)$.
For $k=0$, the same is true if we remove from the above basis the series $E_{2,1,1,t}$ and add instead the series 
$E_{2,\ell,t}$ for $\ell t |N$.
\end{prop}

We now construct a basis of the space of boundary modular symbols similar to the basis of Eisenstein series given in Proposition~\ref{p:miyakestein}.

Let $M$ be a positive integer.
Let $u$ and $v$ be relatively prime integers.  Define $\phi_{k, u, v} \in \Hom(\Delta, \V_k)^{\Gamma_1(M)}$ to be supported on the 
$\Gamma_1(M)$-orbit of $u/v \in \PP^1(\Q)$ and given on that orbit by the formula
\begin{equation} \label{e:phiuvdef}
 \phi_{k, u,v}\left(\gamma\left(\frac{u}{v}\right)\right)(P(z)) = P\left(\gamma\left(\frac{u}{v}\right)\right) \cdot (cu + dv)^k, \qquad \gamma = \mat{a & b \\  c  & d} \in \Gamma_1(M).
\end{equation}
The right side of (\ref{e:phiuvdef}) has the expected meaning when $\gamma(u/v) = \infty$, i.e.\ the value $a_k(au + bv)^k$ if $P(x) = a_k x^k + a_{k-1}x^{k-1} + \cdots.$
One must check that (\ref{e:phiuvdef}) is well-defined, i.e.\ that if $\gamma \in \Gamma_1(M)$ stabilizes $u/v$, then the value of (\ref{e:phiuvdef}) is equal to $P(u/v)v^k$.
It is easy to check that this holds unless $M \mid 4$ and $k$ is odd; since there are no odd characters of conductor 1 or 2, the only problem occurs when $M = 4$ and $k$ is odd;
in this case  one furthermore sees that a problem occurs only when $v \equiv 2 \pmod{4}$; we will not need to define $\phi_{k, u,v}$ in this case.

Note that \begin{equation} \label{e:phiuvsign}
\phi_{k, -u,-v} = (-1)^k \phi_{k, u,v}. \end{equation}
  Furthermore
\begin{equation} \phi_{k, au+bv, cu+dv} = \phi_{k, u,v}, \qquad \mat{a & b \\  c  & d} \in \Gamma_1(M). \label{e:phiuvinv}
\end{equation}
In particular, $\phi_{k, u,v}$ depends only on the $\Gamma_1(M)$-orbit of $u/v$ when $k$ is even, and the same is true up to sign if $k$ is odd.

Now let $\psi$ and $\tau$ be Dirichlet characters of conductors $Q, R$, respectively, with $QR=M$ and $\psi\tau(-1) = (-1)^k$.
Define $ \phi_{k, \psi, \tau} \in \Hom(\Delta, \V_k)^{\Gamma_1(M)}$ by
\begin{equation} \label{e:phipsitau} \phi_{k, \psi, \tau} =
 \sum_{\genfrac{}{}{0pt}{}{x \!\!\!\! \pmod{Q}}{(x,Q)=1}}
 \sum_{\genfrac{}{}{0pt}{}{y \!\!\!\! \pmod{R}}{(y,R)=1}}
\psi^{-1}(x) \tau(y) \phi_{k, x, Qy}.\end{equation}
Equations (\ref{e:phiuvsign}) and (\ref{e:phiuvinv}) together with $\psi\tau(-1) = (-1)^k$ imply that each summand $ \psi^{-1}(x) \tau(y) \phi_{k, x, Qy}$
depends only on the cusp of $\Gamma_1(M)$ determined by $x/Qy$, i.e.\ only on $x \pmod{Q}$ and $y \pmod{R}$. 
 We omit the proof of the following 
proposition, which is a simple computation.

\begin{prop} \label{p:boundsymb}
The vector $\phi_{k, \psi, \tau} \in \Hom(\Delta,\V_k)^{\Gamma_1(M)}$ is an eigenvector for the operators 
\begin{itemize}
\item[(i)]
$T_\ell$ for $\ell \nmid M$ $($resp. $U_\ell$ for $\ell \mid M)$ with  eigenvalues $\psi(\ell) + \tau(\ell)\ell^{k+1}$,
\item[(ii)]  $\langle a \rangle$ for  $a \in (\Z/M\Z)^\ast$  with
eigenvalues $\psi(a)\tau(a)$, and
\item[(iii)]  $\iota$ with eigenvalue $\psi(-1)$.
\end{itemize}
The vectors $\phi_{k,\psi,\tau} \in \Hom(\Delta, \V_k)^{\Gamma_1(M)}$ are generators of their eigenspace for the $T_\ell$ $(\ell \nmid M)$ and 
$\langle a \rangle$. 
\end{prop}

By (\ref{bsymbvg}), the map $\Hom(\Delta,V_k)^{\Gamma_1(M)} \rightarrow \BSymb_{\Gamma_1(M)}(\V_k)$ is an isomorphism unless $k=0$; in this case the kernel is  the line generated by $\phi_{0,1,1}$. For $M=\ell$ a prime, let $\phi_{0,\ell} \in \Hom(\Delta,V_k)^{\Gamma_1(\ell)}$ be defined by formula (\ref{e:phipsitau}) with $Q=1$, $R=\ell$, $\psi$ the trivial Dirichlet character and $\tau$ the {\em non-primitive} principal Dirichlet character of modulus $\ell$ but conductor $1$. Then
one checks that $\phi_{0,\ell}$ is an eigenvector for $T_{\ell'}$ ($\ell' \not = \ell$) with eigenvalue $1+\ell'$,  for $\langle a \rangle$ with eigenvalue $1$ (and is the unique such eigenvector, up to scaling, in $\BSymb_{\Gamma_1(\ell)}(\V_0)$), and of $U_\ell$ with eigenvalue $1$. 
The symbol $\phi_{0, \ell}$ has eigenvalue 1 for $\iota$.

\subsection{Boundary symbols in raised level}
\label{ss:bs}

Now let $N$ be a multiple of $M=QR$.  We will be interested studying the subspace of $\BSymb_{\Gamma_1(N)}(\V_k)$
on which the Hecke operators $T_\ell$ for $\ell \nmid N$ and $\langle a \rangle$ for $a \in (\Z/N\Z)^\ast$ act via the corresponding
Hecke eigenvalues of $E_{k+2, \psi, \tau}$; this subspace will be denoted 
 \[ \BSymb_{\Gamma_1(N)}(\V_k)[E_{k+2, \psi, \tau}]. \]
 
The statement of the following proposition follows directly from Proposition~\ref{p:bsymb}(iii).  However, we provide a separate, more computational proof that
makes explicit the isomorphism (\ref{e:bsymbeis}).

\begin{prop} \label{c:isombsymb}  Let $k \ge 0$ be an integer, and let $\psi$ and $\tau$  be Dirichlet characters of conductors $Q, R$, respectively, such that $E_{k+2, \psi, \tau}$ is a normal Eisenstein series.
Let $M=QR$. The dimension of the space 
 \[  \BSymb_{\Gamma_1(N)}(\V_k)[E_{k+2, \psi, \tau}] \] is equal to the number of divisors of $N/M$.
More precisely,
there exists an isomorphism \begin{equation} \label{e:bsymbeis}
 \BSymb_{\Gamma_1(N)}(\V_k)[E_{k+2, \psi, \tau}] \simeq \EE_{k+2}(\Gamma_1(N))[E_{k+2, \psi, \tau}], 
 \end{equation}
compatible with  the Hecke operators $T_\ell$ for $\ell \nmid N$, $U_\ell$ for $\ell \mid N$, and $\langle a \rangle$ for $a \in (\Z/N\Z)^\ast$.
The same is true for the exceptional Eisenstein series $E_{2, \ell}$, with $M$ replaced by $\ell$.
\end{prop}

\begin{pf}
We make the following simplifying assumption that is sufficient for the purposes of this paper:
\begin{equation} \label{a:sf} N/M \text{ is square-free and relatively prime to } M. \end{equation}
 With this assumption, the space
\begin{equation} \label{e:eiseig}
 \EE_{k+2}(\Gamma_1(N))[E_{k+2, \psi, \tau}] 
 \end{equation}
is semisimple as an algebra for the Hecke operators listed in the Proposition.  To be precise, for each
prime $\ell$ dividing $N/M$, we define Hecke operators providing the $\ell$-ordinary and $\ell$-critical stabilizations
of $E_{k+2, \psi, \tau}$ as follows:
\[ O_\ell := 1 - \tau(\ell)\ell^{k+1}V_\ell, \quad
C_\ell := 1- \psi(\ell)V_\ell.
\]
 For each factorization $N/M = st$ into positive integers, we define
 \begin{equation} \label{e:eisst}
  E_{k+2, \psi, \tau}^{s, t} := (E_{k+2, \psi, \tau})_{|\prod_{\ell \mid s} O_\ell \prod_{\ell \mid t} C_\ell}. 
  \end{equation}
 This form is the $s$-ordinary, $t$-critical eigenvector, i.e.\ the action of the Hecke operator $U_\ell$ for
$\ell \mid N/M$ is given by
\begin{equation}
 (E_{k+2, \psi, \tau}^{s, t})_{| U_\ell} =  E_{k+2, \psi, \tau}^{s, t} \cdot \begin{cases} \psi(\ell) & \ell \mid s \\ \tau(\ell)\ell^{k+1} & \ell \mid t. \end{cases}
\label{e:eissordtcrit}
\end{equation}
The space  (\ref{e:eiseig}) has as a basis the $E_{k+2, \psi, \tau}^{s, t}$. 

It remains to find a basis for $\BSymb_{\Gamma_1(N)}(\V_k)[E_{k+2, \psi, \tau}]$ consisting of eigenvectors with these same eigenvalues.
A natural idea would be to define $\phi_{k+2, \psi, \tau}^{s, t}$ from $\phi_{k+2, \psi, \tau}$  following (\ref{e:eisst}).  However, as we shall see, it is possible for $O_\ell$ to annihilate $\phi_{k+2, \psi, \tau}$.  For this reason we must employ a more explicit approach.

To this end, define $\phi_{k, u, v}^s \in \Hom(\Delta, \V_k)^{\Gamma_1(M) \cap \Gamma_0(s)}$ to be the boundary symbol supported on the $\Gamma_1(M)\cap \Gamma_0(s)$-orbit of $u/v$, and defined on that orbit 
by equation (\ref{e:phiuvdef}), with $\gamma$ restricted to lie in $\Gamma_1(M)\cap \Gamma_0(s)$.
We then define
\begin{equation} \label{e:phipsitaus} \phi_{k, \psi, \tau}^{s} =
 \sum_{\genfrac{}{}{0pt}{}{x \!\!\!\! \pmod{Q}}{(x,Q)=1}}
 \sum_{\genfrac{}{}{0pt}{}{y \!\!\!\! \pmod{R}}{(y,Rs)=1}}
 \psi^{-1}(x) \tau(y) \phi_{k, x, Qy}^s.\end{equation}
Finally, define $\phi_{k, \psi, \tau}^{s, t} \in \Hom(\Delta, \V_k)^{\Gamma_1(M) \cap \Gamma_0(N)}$ by
\[ \phi_{k, \psi, \tau}^{s, t} := (\phi_{k, \psi, \tau}^{s})_{| \prod_{\ell \mid t} C_\ell}. \]

It is not difficult to verify that $\phi_{k, \psi, \tau}^{s, t}$ is $s$-ordinary and $t$-critical, i.e.\ it has eigenvalues for the $U_\ell$, $\ell \mid N/M$, given by (\ref{e:eissordtcrit}).  We must verify that the vectors  $\phi_{k, \psi, \tau}^{s, t}$ are nonzero.
 It is clear from the definition (\ref{e:phipsitaus}) that
$ \phi_{k, \psi, \tau}^{s} \neq 0$. A simple calculation shows that for $\ell \nmid Ms$ and
a cusp $c = x/y$ in lowest form,
we have
\begin{equation} \label{e:vlbound}
 (\phi_{k, \psi, \tau}^s)_{| V_\ell}(c) = \phi_{k, \psi, \tau}^{s}(c) \cdot \begin{cases} \psi^{-1}(\ell)\ell^{-1} & \ell \nmid y \\
\tau^{-1}(\ell)\ell^{-k-1} & \ell \mid  y.
\end{cases} 
\end{equation}
It follows that on the $\Gamma_0(\ell)$-orbit of $0$, we have 
$ (\phi_{k, \psi, \tau}^s)_{| C_\ell}  = (1 - \ell^{-1}) \phi_{k, \psi, \tau}^s \neq 0$;
Therefore, the $\phi_{k, \psi, \tau}^{s, t}$ are nonzero and hence constitute a basis of eigenvectors for 
$\BSymb_{\Gamma_1(N)}(\V_k)[E_{k+2, \psi, \tau}]$.

The result now follows with the isomorphism (\ref{e:bsymbeis}) given by $\phi_{k, \psi, \tau}^{s, t} \mapsto E_{k, \psi, \tau}^{s, t}.$
 This proof is valid for the exceptional Eisenstein series $E_{2, \ell}$ as well, with $\phi_{k, \psi, \tau}$ replaced by $\phi_{0, \ell}$. 
\end{pf}

Equation (\ref{e:vlbound}) implies that for $\ell \nmid Ms$, we have
\[ (\phi_{k, \psi, \tau}^s)_{| O_\ell} = (1 - \tau\psi^{-1}(\ell) \ell^k) \phi_{k, \psi, \tau}^{s \ell}. 
\]
The possible vanishing of the factor on the RHS explains why it was necessary to define $\phi_{k, \psi, \tau}^{s}$ with the explicit
formula (\ref{e:phipsitaus}), rather than simply as $(\phi_{k, \psi, \tau})_{\prod_{\ell \mid s}O_\ell}$, 
and motivates the following definition.

 \begin{definition} \label{c:bad}
 Let $f$ be a new Eisenstein series of level $\Gamma_1(M)$, and $\ell$ be a prime not dividing $M$.
 We say that $\ell$ is {\it bad} for $f$ if either $f$ is exceptional, or $f$ is the normal Eisenstein series $E_{k+2, \psi, \tau}$ with $k=0$ and $\tau(\ell) = \psi(\ell)$.
 \end{definition}
 
   If we assume that the primes dividing $N/M$ are not bad for
 $E_{k+2, \psi, \tau}$, then we arrive at the following simpler version of Proposition~\ref{c:isombsymb}.

\begin{prop}  Let $E_{k+2, \psi, \tau}$ be a normal Eisenstein series of level $M=QR$ dividing $N$.
Suppose that $N/M$ is squarefree, relatively prime to $M$, and such that all primes $\ell$ dividing $N/M$ are not bad for  $E_{k+2, \psi, \tau}$.
Then the symbols $(\phi_{k, \psi, \tau})|_{V_t}$ for positive integers $t$ dividing $N/M$ provide a basis for 
$\BSymb_{\Gamma_1(N)}(\V_k)[E_{k+2, \psi, \tau}] $,
and the
map $(\phi_{k, \psi, \tau})_{|V_t}  \mapsto (E_{k+2, \psi, \tau})_{|V_t}$ 
yields an isomorphism \[ \BSymb_{\Gamma_1(N)}(\V_k)[E_{k+2, \psi, \tau}] \simeq \EE_{k+2}(\Gamma_1(N))[E_{k+2, \psi, \tau}] \]
that is compatible with  the Hecke operators $T_\ell$ for $\ell \nmid N$, $U_\ell$ for $\ell \mid N$, and $\langle a \rangle$ for $a \in (\Z/N\Z)^\ast$.
\end{prop}


We conclude this section by recording a corollary of these computations that gives a basis for the space of ordinary boundary symbols.
Fix an integer $p$ not dividing $N$ and write $\Gamma = \Gamma_1(N) \cap \Gamma_0(p)$.
Write $\BSymb_{\Gamma}(\V_k(L))^{\ord}$ for the subspace generated by the generalized eigenvectors for $U_p$ with eigenvalue a root of unity (and similarly for $ \EE_{k+2}(\Gamma)^{\ord}$).   Note that if $L$ is a finite extension of $\Q_p$, then this definition is unchanged if ``a root of unity" is replaced by ``of $p$-adic order 0." 

\begin{cor} \label{c:ordphi}  Let $L$ be a field of characteristic 0 containing the $\varphi(N)$th roots of unity.  Let $E_{k+2, \psi, \tau}$ be a normal Eisenstein series of level $M=QR$ dividing $N$.  
Suppose that $N/M$ is squarefree, relatively prime to $M$, and such that all primes $\ell$ dividing $N/M$ are not bad for  $E_{k+2, \psi, \tau}$.
Let $p$ be a prime not dividing $N$.
Then the symbols $\phi_{k, \psi, \tau, t}^p:= (\phi_{k, \psi, \tau}^p)|_{V_t}$ for positive integers $t$ dividing $N/M$ provide a basis for 
 $\BSymb_{\Gamma}(\V_k(L))^{\ord}[E_{k+2, \psi, \tau}]$, and the map
 \[ \phi_{k, \psi, \tau,t}^p \mapsto E_{k+2, \psi, \tau,t}^{\ord}:= (E_{k+2, \psi, \tau}^{\ord})_{|V_t} \]
yields a Hecke-equivariant isomorphism
\[ \BSymb_{\Gamma}(\V_k)^{\ord}[E_{k+2, \psi, \tau}] \simeq \EE_{k+2}(\Gamma)^{\ord}[E_{k+2, \psi, \tau}].\]
\end{cor}

For future reference, it is useful to note that the vectors $\phi_{k, \psi, \tau,t}^p$ appearing in Corollary~\ref{c:ordphi}
satisfy
$\phi_{k, \psi, \tau,t}(c) = 0$ if $c \not \in \Q \cap \Z_p$, and $\phi_{k, \psi, \tau,t}^p(c)$ is a multiple
of the linear form $P \mapsto P(c)$ if $c \in \Q \cap \Z_p$.

\section{Eisenstein series and their non-boundary modular symbols} \label{s:eismodsymb}

\subsection{Notations} \label{s:notation}

We now fix notations that will stay in force for the remainder of this paper.
We fix an integer $N$ and a prime $p$ not dividing $N$. We set
$$\Gamma := \Gamma_1(N) \cap \Gamma_0(p).$$
We define the set of cusps \[ C := \Gamma_0(N)\infty \cup \Gamma_0(N)0,\] 
 the subset of $\PP^1(\Q)$ containing $\infty$ and all rationals
$\frac{a}{m}$  in lowest terms with  either $N \mid m$ or $(N, m) = 1$.
Note that $\Gamma_0(N)$, and hence its subgroup $\Gamma$, stabilizes $C$.  Observe that the matrices $\mat{ 1 & 0 \\ 0 & d}$ also stabilize $C$ whenever $d \in \Z$ is prime to $N$.  Let $S$ be the submonoid 
of $\Gl_2(\Q)$ generated by those matrices and $\Gamma_0(Np)$, so $S$ stabilizes $C$ as well.

We shall consider in this paper spaces of partial modular symbols $\Symb_{\Gamma,C}(V)$ (and $\BSymb_{\Gamma,C}(V)$)
for various right $S$-modules $V$. These spaces have an  action of the Hecke operators $T_\ell$ (for $\ell \nmid Np$), $U_p$, and $\langle a \rangle$ for $a \in (\Z/N\Z)^\ast$. We denote by $\HH \subset \HH(\GL_2^+(\Q), \Gamma)$ 
the subalgebra generated by these operators. 
Note that there is also an action of the involution $\iota$ on the space $\Symb_{\Gamma, C}(V)$.

We will also have to consider spaces of full modular symbols $\Symb_\Gamma(V)$ and $\BSymb_\Gamma(V)$. They also have an action of $\HH$ and $\iota$. When $V$ has an action of the larger monoid $S'$ generated by $\Gamma_0(N)$ and all the matrices $\mat{1 & 0 \\ 0 & d}$ for all $d$ in $\Z_{(p)}$, $d \neq 0$ (in practice this will always be the case),
then we also have an action of operators $U_\ell$ (for $\ell \mid  N$) on $\Symb_\Gamma(V)$, and maps
\[ V_t: \Symb_{\Gamma_1(M) \cap \Gamma_0(p)}(V) \rightarrow \Symb_\Gamma(V) \] for integers $M$ and $t$
such that $Mt \mid N$. These operators and maps will be useful in \S\ref{s:families}.

\subsection{$C$-cuspidal modular forms and their symbols}

We now turn to the more interesting modular symbols associated to modular forms,
in contrast to the boundary symbols studied in \S\ref{s:bound}.  These symbols are
 defined by
integration on the upper half-plane and are related to special values of $L$-functions.
Recall that if $f$ is a modular form for $\Gamma$, then $f$ has a $q$-expansion
$$f(z) = \sum_{n=0}^{\infty}  a_n q^n,\ \ \ q=e^{2 i \pi z}.$$
For any Dirichlet character $\chi$, we define the 
{\em $L$-function of $f$ twisted by $\chi$} as the Dirichlet series
 \begin{equation}
 \label{e:lfchi}
  L(f,\chi,s)= \sum_{n=1}^\infty a_n \chi(n) n^{-s}. 
  \end{equation}
  This converges for $\Re(s)$ large  enough and admits a meromorphic continuation to the complex plane. Note that we do not include $a_0$ in the sum.  More generally, for any positive integer $d$, define
 \[ L(f, \chi, d, s) = d^{1-s} \sum_{n=1}^\infty a_{nd} \chi(n) n^{-s}, \]
 so $(\ref{e:lfchi})$ is the case $d=1$.

Let $f \in M_{k+2}(\Gamma,\C)$ be a modular form of weight $k+2$ for $\Gamma$.
We say that $f$ is {\it $C$-cuspidal} if it vanishes at all cusps in $C = \Gamma_0(N) \infty \cup \Gamma_0(N)0$.  We henceforth assume that $f$ is $C$-cuspidal.

Let $m$ be an integer relatively prime to $N$ and $\chi$ a Dirichlet character of conductor dividing $m$, say $m = \cond(\chi) \cdot d$. We have the well-known and elementary formula
\begin{eqnarray} \label{Lfchis} 
\frac{\Gamma(s)}{(2\pi)^s}L(f,\chi,d, s) = \frac{1}{G(\overline{\chi})} \sum_{\ \ a\!\!\!\!\! \pmod{m}}  \overline{\chi}(a) \int_0^{\infty} f(iy + a/m) y^{s-1} dy,
\end{eqnarray}
where the integral converges because $f$ vanishes at the cusps $a/m$ and $\infty$. 
Here $G(\overline{\chi})$ denotes the Gauss sum. 

\begin{definition} \label{d:defpartial} For $f$ a $C$-cuspidal modular form as above, we define \[ \phi_f \in \Hom(\Delta_C^0,\V_k(\C)) \] by setting \[ \phi_f(\{\infty\}-\{a/m\})(P) = \int_{a/m}^{i\infty} f(z) P(z) dz \] for all $P \in \P_k(\C)$.
\end{definition}

The following proposition is a standard computation.

\begin{prop} We have $\phi_f \in \Symb_{\Gamma,C}(\V_k(\C))$. Moreover
the map $f \mapsto \phi_f$ is $\HH$-equivariant.\end{prop}

Recall that $\Symb_{\Gamma,C}(\V_k(\C))$ is endowed with an action of the involution $\iota$.
We set $$\phi_f^+ = \frac{\phi_f + (\phi_f)_{|\iota}}{2}\quad \text{ and } \quad
\phi_f^- = \frac{\phi_f - (\phi_f)_{|\iota}}{2},$$ so that 
\[ \phi_f = \phi_f^+ + \phi_f^- \quad \text{ and } \quad
(\phi_f^\pm)_{|\iota} = \pm \phi_f^\pm. \]

In view of the definition of $\phi_f$, equation (\ref{Lfchis}) implies immediately that for any integer $j$, $0 \leq j \leq k$, and all Dirichlet characters $\chi$ of conductor  $m/d$:
\begin{equation} \label{formulaLf1}  \frac{j! G(\overline \chi)}{(-2  \pi i)^{j+1}}L(f,\chi,d, j+1) = \!\!\!\!\!\!\!\!  \sum_{\ \ a\!\!\!\!\! \pmod{m}} \!\!\!\!\! \overline \chi(a) \phi_f(\{\infty\}-\{a/m\})( (z-a/m)^j).
 \end{equation}
Noting that the right hand side is multiplied by $\chi(-1)(-1)^j $ if each occurrence of $\phi_f$ in it is replaced by $(\phi_f)_{|\iota}$, we see that for every choice of sign $\varepsilon = \pm 1$, we have
\begin{equation} \label{formulaLf2}  \frac{j! G(\overline \chi)}{(-2 \pi i)^{j+1}}L(f,\chi,d, j+1) = \!\!\!\!\!\!\!\! \sum_{\ \ a\!\!\!\!\! \pmod{m}} \!\!\!\!\! \overline \chi(a) \phi_f^\varepsilon(\{\infty\}-\{a/m\})( (z-a/m)^j)\end{equation}
  for any integer $j$, $0 \leq j \leq k$, and all Dirichlet characters $\chi$ of conductor  $m/d$ such that $\chi(-1)(-1)^j = \varepsilon.$

From (\ref{formulaLf2}) and the linear independence of characters, we deduce
\begin{prop} \label{propcriterionsign}  Fix a sign $\varepsilon=\pm 1$. 
\begin{itemize}
\item[(i)] We have $\phi_f = \phi_f^\varepsilon$ $($i.e.\ $\phi_f^{-\varepsilon}=0)$ if and only if for all integers $j$, $0 \leq j \leq k$, and all Dirichlet characters $\chi$ of conductor 
prime to $N$ such that $\chi(-1)(-1)^{j} = - \varepsilon$ and all positive integers $d$, we have $L(f,\chi,d, j+1)=0$.
\item[(ii)] Let $\overline{\Q}$ be the subfield of algebraic numbers in $\C$. We have $\phi_f^\varepsilon \in \Symb_{\Gamma,C}(\V_k(\overline{\Q}))$ if and only if for all integers $j$, $0 \leq j \leq k$,  all Dirichlet characters $\chi$ of conductor prime to $N$ such that $\chi(-1)(-1)^{j} =  \varepsilon$, and all positive integers $d$, we have $L(f,\chi,d,j+1) / \pi^{j+1} \in \overline{\Q}$. 
\end{itemize}

\end{prop}

\subsection{$C$-cuspidal Eisenstein series}  \label{s:ccuspidal}
We now apply the observations of the previous section to the case of $C$-cuspidal Eisenstein series.

\begin{lemma} \label{Leisensteinlemma} 
 We have the identity
$$ L(E_{k+2,\tau,\psi,t}, \chi, d, s) = \chi(t) t^{-s} \tau(d)d^{1-s} \sum_{u \mid d} u^{k+1} \psi\tau^{-1}(u) L(\tau \chi,s) L_{d/u}(\psi \chi,s-k-1),$$
where $L_{d/u}(\psi \chi,s-k-1)$ denotes the Dirichlet $L$-function with the Euler factors at primes dividing $d/u$ removed.
\end{lemma}
\begin{pf} 
We prove the simpler formula
\[  L(E_{k+2,\tau,\psi,t}, \chi, s) = \chi(t) t^{-s}  L(\tau \chi,s) L(\psi \chi,s-k-1)
\]
for the case $d=1$, leaving the general case to the reader.

It suffices to prove the result for  $\Re(s)$ large enough.
If $a=(a_n)_{n \geq 1}$ is a sequence of complex numbers, the $L$-function of $a$ is by definition the function $L(a,s)=\sum_{n \geq 1} a_n/n^s$, which converges on some ``half-plane" 
$\Re(s) > \rho$ for some $\rho \in \R \cup \{-\infty,+\infty\}$. If $(a_n)$ and $(b_n)$ are two sequences, one defines a commutative and associative convolution product  $(a_n) \ast (b_n)$ as the sequence $(c_n)$ given by $c_n=\sum_{d|n} a_d b_{n/d}$. It is then elementary and well-known that 
$L(c,s)=L(a,s)L(b,s)$ on the domains of convergence of those series. Note that if
$a,b$ are two sequences, and  $c$ is strictly multiplicative (that is, $c_{nm}=c_nc_m$ for any $n,m \geq 1$), then $(ac) \ast (bc) = (a \ast b)c$.

The sequence of coefficients $(a_n)$ of $E_{k+2,\tau,\psi,t}$ is the convolution product $(b_n)\ast(c_n)\ast(d_n)$ where $b,c,d$ are the sequences defined by the equalities (for $n \geq 1$)
 $b_n=\psi(n) n^{k+1}$, $c_n=\tau(n)$, and 
 \[ d_n=\begin{cases} 1 & \text{if } n=t, \\ 0 & \text{if } n \neq t. \end{cases} \]
 The $L$-function $L(E_{k+2,\tau,\psi,t}, \chi, s)$ is by definition the $L$-function of 
 the sequence $a \chi$, which is $(b \chi) \ast (c \chi) \ast (d \chi)$ since $\chi$ is 
 strictly multiplicative.
 But we have by definition $L(b \chi,s)=L(\psi \chi ,s-k-1)$, $L(c\chi,s)=L(\tau \chi,s)$ and $L(d \chi,s)=\chi(t) t^{-s}$.
 The result follows.
 \end{pf} 
 
Now consider a linear combination $f= \sum_{t,\ QRt|N} c_t E_{k+2,\tau,\psi,t}$ where the $c_t$ are algebraic numbers.
If $f$ is $C$-cuspidal 
 then we can attach to $f$ a partial modular symbol $\phi_f$ over $C$.
 
 \begin{prop} \label{symbeisgood}
 Let $\epsilon = -\tau(-1)$. Then $(\phi_f)_{|\iota} = \epsilon \phi_f$ unless
 $\tau=1$ or $\psi=1$; in  these cases, $\phi_f^{-\epsilon}$ is a multiple of the boundary symbol $\phi_{k, \tau, \psi}$ defined in Section~\ref{s:bound}.
  In all cases, $\phi_f \in \Symb_{\Gamma,C}(\V_k(K))$ for some number field $K$.\end{prop}
\begin{pf} According to Proposition~\ref{propcriterionsign}, for the first assertion it suffices to prove that for all integers $j$ such that $0 \leq j \leq k$, all Dirichlet characters $\chi$ of conductor prime to $N$ such that $\chi(-1)(-1)^{j}= \tau(-1),$ and all positive integers $d$, we have
$L(f,\chi,d, j+1)=0$. But by Lemma~\ref{Leisensteinlemma}, the value $L(f,\chi,d, j+1)$ is
a linear combination of terms of the form 
$$ L(\tau \chi,j+1) L_u(\psi \chi,j-k)$$
for divisors $u \mid d$.
By the functional equation for Dirichlet $L$-functions, $L(\psi \chi,j-k)$ vanishes if $j \leq k$ and
 $\psi\chi(-1) (-1)^{j-k} = 1$, except in the case $\psi=\chi=1, j=k$. The condition  $\psi\chi(-1) (-1)^{j-k} = 1$ is equivalent to \[ \chi(-1)(-1)^{j}
 = \psi(-1) (-1)^k = \tau(-1). \]
We therefore obtain our desired vanishing $L(f,\chi,d, j+1)=0$ except in the case $\psi=\chi=1, j=k$ and the case where $L(\tau\chi, j+1)$ has a pole, namely
$\tau = \chi = 1, j=0$.  We leave to the reader the remaining verification from what we have already proved that in the cases $\tau=1$ or $\psi=1$, the symbol
 $\phi_f^{-\epsilon}$ is a multiple of the boundary symbol $\phi_{k, \tau, \psi}$; we will not use this result in the sequel.

To prove that $\phi_f \in \Symb_{\Gamma,\C}(\V_k(\overline{\Q}))$ (from which it follows easily by the type finiteness
of $\Delta_{C}^0$ as a $\Z[\Gamma]$-module that $\phi_f \in \Symb_{\Gamma,\C}(\V_k(K))$ for some number field $K$),
it suffices to check that $L(f,\chi,d,j+1)/(2 \pi i)^{j+1}$ are algebraic for $0 \leq j \leq k$, all Dirichlet characters $\chi$ of conductor prime to $N$ and all positive integers $d$.
But  ${L(\tau\chi,j+1)}/{(2  \pi i)^{j+1}}$ and  $L(\psi \chi,j-k)$ are algebraic by  well-known properties of Dirichlet $L$-functions.  In fact, one can show that $\phi_f$ is defined over the field $\Q(\tau, \psi)$ defined by the values of the characters $\tau, \psi$.
 \end{pf}
 
\section{Overconvergent partial modular symbols}

We fix embeddings $\overline \Q \subset \overline \Q_p$ and $\overline \Q \subset \C$. We denote by $\ord_p\colon \overline \Q_p^\ast \rightarrow \Q$ the $p$-adic valuation, normalized such that $\ord_p(p)=1$. We define  weight space $\WW$ as the rigid analytic space over $\Q_p$ such that $\WW(L) = \Hom_{\cont}(\Z_p^\ast,L^\ast)$ for every $\Q_p$-Banach algebra $L$.
We see $\Z$ as a (Zariski-dense) subset of $\WW(\Q_p)$ by identifying $k \in \Z$ with the character $z \mapsto z^k$.

\subsection{Modules of distributions}

We recall definitions and notations from \cite[\S3]{Bcrit}, which are generalizations of earlier definitions of Stevens (\cite{stevens}).

If $L$ is any $\Q_p$-Banach algebra, and $r>0$, we write $\AA[r](L)$ for the $L$-Banach module
of functions $f\colon \Z_p \rightarrow L$ that are locally analytic of radius $r$ (that is,
such that for any $e \in \Z_p$ there exists a formal series $\sum_{n=0}^\infty a_n (z-e)^n$ with radius of convergence at least $r$ and converging to $f$ on the closed ball of center $e$ and radius $r$ in $\Z_p$). We write $\DD[r](\Q_p)$ for the Banach dual of $\AA[r](\Q_p)$,  and define a Banach $L$-module \[ \DD[r](L) := \DD[r](\Q_p) \hotimes_{\Q_p} L.\]  For $r \ge 0$, we write $\DD^\dagger[r](L)$ for the projective limit of $\DD[r'](L)$ as $r' \rightarrow r^+$, which is a (Frechet) $L$-module. 

Let $\kappa \in \WW(L)$, i.e.\ $\kappa$ is a continuous character $\Z_p^\ast \rightarrow L^\ast$.  There exists $r(\kappa) \in \R_{>0} \cup \{ \infty\}$ such that for all $0 < r \leq r(\kappa)$,
we can provide $\AA[r](L)$
with an action of \[ S_0(p)=  \left\{ \gamma = \left( \begin{array}{cc}      a  & b \\
      c & d \\  \end{array} \right) \in M_2(\Z_{(p)}),\ p \nmid a, \ p \mid c, \ ad-bc \neq 0 \right\},\] which is said to be {\it of weight $\kappa$}, by the rule
\begin{eqnarray} \label{actakappa}
(\gamma \cdot_\kappa f)(z) = \kappa(a-cz)    f\left(\frac{dz-b}{a-cz}\right).
\end{eqnarray}
These actions are compatible for various $r$. By duality they induce actions on $\DD[r](L)$ for $0<r \leq r(\kappa)$, hence on $\DD^\dag[r](L)$ for $0 \leq r < r(\kappa)$.
We write $\DD_\kappa[r](L)$ and $\DD^\dag_\kappa[r](L)$ for the modules  $\DD[r](L)$
and $\DD^\dag[r](L)$ together with the weight $\kappa$ action of $S_0(p)$. (See 
\cite[\S3]{Bcrit} for more details.) For $\kappa = k \in \Z$ we can take $r(k) = \infty$.

Of these spaces, the one we shall use most is $\DD^\dag[0]$. To lighten notation, 
we write $\DD^\dag$ for $\DD^\dag[0]$.

Recall (\cite{stevens}) that when $k$ is a non-negative integer we have an exact sequence of $S_0(p)$-modules
\begin{equation} 
\label{fundexsymb}
\xymatrix{0 \ar[r] & \DD^\dag_{-2-k}[1](L)(k+1)\ar[r]^{\ \ \ \ \ \ \Theta_k} & \DD_k^\dag[1](L)\ar[r]^{\rho_k} &  \V_k(L) \ar[r] & 0}.
\end{equation} 
The $(k+1)$ after the first term means that the action of $s \in S_0(p)$ is multiplied by $(\det s)^{k+1}$.
The map $\Theta_k$ is the $(k+1)$-th iteration of the derivation of distribution, i.e.\  $(\Theta_k \mu)(\sigma) = \mu\!\left(\frac{d^{k+1}\sigma}{dz^{k+1}}\right)$.
The map $\rho_k$ is the dual of the obvious inclusion map 
$\P_k(L) \rightarrow \AA_k[1](L)$. 

\subsection{Spaces of overconvergent partial modular symbols}

The monoid we called $S$ in \S\ref{s:notation} is contained in $S_0(p)$. Therefore, as explained {\it loc.\ cit.},  the modules of partial modular symbols $\Symb_{\Gamma,C}(\DD^\dag_\kappa[r](L))$ (or the same without the $\dag$) have a natural action of the Hecke algebra $\HH$ and the involution $\iota$.

If $V$ is any $\Q_p$-vector space on which an operator denoted $U_p$ acts, and $\nu \in \R$,
we define the {\it slope $< \nu$ part of $V$}, denoted  $V^{< \nu}$, as the sum of the indecomposable $U_p$-submodules of $V$ on which $U_p$ acts (after extending scalars) by eigenvalues of $p$-adic valuation less than $\nu$. 

\begin{prop} \label{symbspace} Assume that $L$ is a finite extension of $\Q_p$, and that $\kappa$ is an integer $k \ge 0$. 
\begin{itemize} 
\item[(i)] The Banach space $\Symb_{\Gamma,C} (\DD_k[r](L))$ is orthonormalizable, and $U_p$ acts on it as a compact operator. For any $\nu \in \R$, and $0 < r' \le r \leq p$, the restriction map
induces  an isomorphism
\[  \Symb_{\Gamma,C} (\DD_k[r'](L))^{<\nu} 
 \longrightarrow \Symb_{\Gamma,C} (\DD_k[r](L))^{<\nu}, \]
and these spaces are finite-dimensional.
Moreover, for $ 0 \le r < p$, the space $\Symb_{\Gamma,C} (\DD_k^\dag[r](L))^{<\nu}$
is also isomorphic (by the restriction map) to the spaces $\Symb_{\Gamma,C} (\DD_k[r'](L))^{<\nu}$.
\item[(ii)] The exact sequence $(\ref{fundexsymb})$ gives by functoriality a sequence of $\HH$-modules
\begin{equation} \label{fundexseq}
\xymatrix{0 \ar[r] &  \Symb_{\Gamma,C}(\DD^\dag_{-2-k}[1](L))(k+1)\ar[r]^{\ \ \ \ \ \ \Theta_k} & \Symb_{\Gamma,C} (\DD^\dag[1]_k(L))\ar[r]^{\rho_k} & \Symb_{\Gamma,C}(\V_k(L)) \ar[r] & 0}
\end{equation}
 that is still exact,
where the $(k+1)$ in the first term means that we multiply the action of $[\Gamma s \Gamma]$ by $(\det s)^{k+1}$.
\item[(iii)] $($Stevens' control theorem$)$ For any $r$ such that $0 \leq r < p$, the map
 $$  \xymatrix{\Symb_{\Gamma,C} (\DD^\dag_k[r](L))^{<k+1} \ar[r]^{\ \ \rho_k} & \Symb_{\Gamma,C}(\V_k(L))^{<k+1}}$$
  given by the composition of $\rho_k$ from (ii) and the isomorphism from (i)
  is an isomorphism.
\end{itemize}
 \end{prop}
\begin{pf}
These results are proved exactly as their counterparts for full modular symbols are proved in \cite{stevenspollack2}, using for (ii) the cohomological interpretation given in Proposition~\ref{propisomcohom}.
\end{pf}
\subsection{The $p$-adic $L$-function of an overconvergent
partial modular symbol}                  

We shall denote by $\RR$ the $\Q_p$-Frechet space $\anneau(\WW)$ of the rigid analytic space $\WW$.
Given a distribution $\mu \in \DD^\dagger(\Q_p)$, we define its {\it Mellin transform} $\Mel(\mu) \in \RR$ by $\Mel(\mu)(\sigma) = \int_{\Z_p^\ast} \sigma d\mu$. 
This makes sense since a character $\sigma \in \WW(\C_p)$ is always a locally analytic function of some radius. 

The Mellin transform realizes a bijection (even a Frechet isomorphism) between the subspace of $\DD^\dagger(\Q_p)$ of distributions with support in $\Z_p^\ast$ and $\RR$. But note that distributions with support in $p \Z_p$ have zero Mellin transform.

More generally, for any $\Q_p$-Banach algebra $L$ we have a Mellin transform map $\Mel_L: \DD^\dagger(L) \rightarrow \RR \hotimes L$ (cf. \cite{Bcrit}).

\begin{definition} If $L$ is a $\Q_p$-Banach
algebra, we define a map \[ L_p :  \Symb_{\Gamma,C} (\DD_\kappa^\dagger(L)) \rightarrow
\RR \hotimes L\] by \[ L_p(\Phi) = \Mel_L( \Phi(\{\infty\} - \{0\})). \]
\end{definition}

\subsection{The eigencurves of partial modular symbols}

\begin{prop} \label{symbmodule} Let $X=\sp A$ be an open affinoid subset of $\WW$, and let $K \in \WW(A)$ be the canonical character, i.e.\ the $A$-point of $\WW$ corresponding to the inclusion $X=\sp A \subset \WW$. We write $r(X)$ for $r(K)$.
\begin{itemize}
\item[(i)] For $0<r \leq r(X)$, the modules $\Symb_{\Gamma,C}( \DD_K[r] (A))$ are orthonormalizable Banach $R$-modules and $U_p$ acts compactly on them.  These modules for different $r$ are {\it linked} in the sense of \cite[\S5]{Buz}. The formation
of these modules commutes with restriction to an open affinoid subdomain: if $\sp A'=X' \hookrightarrow X$ is an open affinoid subdomain, then 
$$\Symb_{\Gamma,C}( \DD_K[r] (A)) \hotimes_A A' =  \Symb_{\Gamma,C}( \DD_K[r] (A'))$$
\item[(ii)]  If $\nu \in \R$ and $X$ is small
enough (with respect to $\nu$, but independently of $r$), then one can naturally define a sub-module $\Symb_{\Gamma,C}( \DD_K[r](A))^{< \nu}$ of $\Symb_{\Gamma,C}(\DD_K[r](A))$, which is a finite projective module over $A$. The restriction maps between these modules for various $r$ are isomorphisms, and we can therefore  identify these modules
with a submodule  $\Symb_{\Gamma,C}( \DD_K[r]^\dag (A))^{< \nu}$ of $\Symb_{\Gamma,C}( \DD_K^\dag(A))$
for $0 \le r < r(X)$.
\item[(iii)] Let $\kappa \in X$ be a point with field of definition $L(\kappa)$, a finite extension of $\Q_p$. So $\kappa \in X(L(\kappa)) \subset \WW(L(\kappa))$.
The natural specialization morphism of $\HH$-modules
 \begin{eqnarray} \label{specmap} \Symb_{\Gamma,C}( \DD_K^\dagger(A))^{< \nu}  \otimes_A L(\kappa) 
 \longrightarrow \Symb_{\Gamma,C}( \DD_\kappa^\dagger(L(\kappa)))^{< \nu} \end{eqnarray}
is injective,  has corank at most one, and is surjective except possibly if $\kappa$ is the trivial character.
\end{itemize}
\end{prop}

\begin{pf} The results of (i) and (iii) are proved exactly as those for full modular symbols:
see \cite[\S3]{Bcrit} or \cite{Bcourse}. The results of (ii) follow from those of (i)
and Coleman--Buzzard's Riesz theory (\cite[\S5]{Buz}, \cite[\S3]{Bcrit} and \cite{Bcourse}).
\end{pf}

We can use the modules $\Symb_{\Gamma,C}( \DD_K[r](A))$ to construct two eigencurves of partial modular symbols.
To be precise, we split 
$\Symb_{\Gamma,C}( \DD_K[r](A))$ into two eigenspaces $\Symb_{\Gamma,C}( \DD_K[r](L))^{\pm}$, with eigenvalue $\pm 1$ for the involution $\iota$. 
We fix a choice of sign $\pm$ and apply Buzzard's eigenvariety machine \cite[Construction 5.7]{Buz}. For any admissible affinoid subset  $X=\sp A \subset \WW$,  we write $M_X$ for the module 
$\Symb_{\Gamma,C}( \DD_K[r] (A))^{\pm}$ where we have chosen an $r \leq r(X)$. For Buzzard's  commutative Hecke algebra ${\mathbf T}$, we take our algebra $\HH \subset \HH(\GL_2^+(\Q), \Gamma)$,
and for Buzzard's operator $\phi \in {\mathbf T}$, we take the operator $U_p$. By Prop.~\ref{symbmodule}(i), $M_X$ satisfies Buzzard property (Pr) (see \cite[\S2]{Buz})and the operator $\phi=U_p$ acts compactly on it. Moreover,
 $M_X \hotimes_{A} A'$ is linked (in the sense of \cite[section 5]{Buz}) to $M_{X'}$ if $X'=\sp A'$ is a subdomain of $A$,
 as can be proved exactly as the analog statement for full modular symbols in \cite{Bcourse}. Therefore we can apply Buzzard's eigenvariety machine
to define {\em the eigencurve $\CC_{\Gamma,C}^\pm$ of partial modular symbols} for $C$ and $\Gamma$ of sign $\pm$.  It is independent of the choices we have made (of an $r \leq r(X)$ for every $X$).

The eigencurve $\CC_{\Gamma,C}^\pm$ enjoys the following usual properties of an eigencurve. It is a reduced rigid analytic space of equidimension $1$  over $\Q_p$,
endowed with the following additional structure:
\begin{itemize} \item[(a)] a locally finite flat {\it weight map} $w\colon \CC_{\Gamma,C}^\pm \rightarrow \WW$; 
\item[(b)] a morphism of rings  $\HH \rightarrow \anneau(\CC_{\Gamma,C}^\pm)$ that sends $U_p$ to an invertible (i.e.\ non-vanishing) function.
\end{itemize}
We do not give a name to the morphism (b): an element $t \in \HH$ defines an analytic function on $\CC_{\Gamma,C}^\pm$ that we denote by the same symbol $t$.
If a point $x \in \CC_{\Gamma,C}^\pm$ defined over a finite extension $L$ of $\Q_p$ maps to a weight $\kappa \in \WW(L)$, then $x$ corresponds
to a non-zero $\HH$-eigenspace in $\Symb_{\Gamma,C}(\DD_\kappa^\dag(L))^\pm$,
the eigenvalue of an element $t \in \HH$ being the value $t(x)$ of the function $t$ on the eigencurve at $x$.   Conversely, any eigenspace $\Symb_{\Gamma,C}(\DD_\kappa^\dag(L))^\pm$ on which the $U_p$-eigenvalue is nonzero (except perhaps one in the case $\kappa=0 $) corresponds to
a point $x$ in the eigencurve defined over $L$ and lying over $\kappa$; this is the content of Prop~\ref{symbmodule}(iii) (completed with Prop.~\ref{symbspace}(i)). 

By Prop.~\ref{symbspace}(iii),  if a point $x$ as above lies over a weight $\kappa$ that is a non-negative integer $k$, and satisfies 
$\ord_p(U_p(x)) < k+1$, then $x$ corresponds in fact to an eigenspace in the
space of {\it classical} partial modular symbols $\Symb_{\Gamma,C}(\V_k(L))^\pm$.
Such classical points are Zariski-dense in $\CC_{\Gamma, C}^\pm$.

One can use Chenevier's comparison theorem to establish precise relations between
the $\CC_{\Gamma, C}^\pm$ for various choices of sets of cusps $C$ and sign $\pm$,
and also between those curves and the standard Buzzard--Coleman--Mazur eigencurve.
We will content ourselves with the following result, where $\CC_{\Gamma,\BCM}$ denotes the traditional  Buzzard--Coleman--Mazur eigencurve constructed from overconvergent modular forms, and $\CC_{\Gamma,\BCM}^\cusp$ denotes the  cuspidal  Buzzard-Coleman-Mazur eigencurve (i.e.\ the part of the eigencurve parameterizing cuspidal overconvergent modular forms), defined using the same Hecke algebra $\HH$.

\begin{prop} \label{compeigencurveBCM} There are  unique closed immersions
(compatible with the maps over $\WW$,
and the map from $\HH$) $$\CC_{\Gamma,\BCM}^\cusp \hookrightarrow \CC_{\Gamma, C}^\pm \hookrightarrow  \CC_{\Gamma,\BCM}.$$
Moreover, if $x$ is a point of $\CC_{\Gamma,\BCM}^\cusp(L)$ for some finite extension  $L/\Q_p$, 
with $w(x)=k \in \Z$ and $v_p(U_p(x)) > 0$, we have 
$$\dim S_{k+2}^\dag(\Gamma,L)_{(x)} = \dim \Symb_{\Gamma,C}^{\pm}(\D_{k}^\dag(L))_{(x)}.$$
\end{prop}
\begin{pf} The first statement is an application of Chenevier's theorem on comparison of eigenvarieties, given the comparison of their {\it classical structure} that we proved in Corollary~\ref{corparsymb}. We do not give the details of the proof as they are  similar to the proof of the special case of full modular symbols given in \cite[Theorem 3.27]{Bcrit}.

For the second, we observe from the same theorem that
$$\dim S_{k+2}^\dag(\Gamma,L)_{(x)} \leq \dim \Symb_{\Gamma,C}^{\pm}(\D_{k}^\dag(L))_{(x)}
\leq \dim M_{k+2}^\dag(\Gamma,L)_{(x)}.$$  Since $x$ is non-ordinary, it does not lie on the Eisenstein components of $\CC_{\Gamma,\BCM}$, so
we have $ S_{k+2}^\dag(\Gamma,L)_{(x)} = M_{k+2}^\dag(\Gamma,L)_{(x)}$ and the result follows.
\end{pf}

\begin{cor} \label{cordimxy} Let $x \in \CC_{\Gamma,\BCM}^\cusp(L)$, with $w(x)=k \in \Z$ and $v_p(U_p(x)) > 0$. Assume further that $x$ is a smooth point of $\CC_{\Gamma,\BCM}$.
Then for any  $y \in  \CC_{\Gamma,\BCM}^\cusp(L)$ sufficiently close to $x$ and such that $w(y) = k' \in \N$ and $v_p(U_p(y)) < k'+1$ $($that is, $y$ satisfies the hypothesis of Coleman's control theorem and is therefore classical$)$, we have
$$ \dim \Symb_{\Gamma,C}^{\pm}(\D_{k}^\dag(L))[x] = \dim {S_{k+2}(\Gamma,L)[y]}.$$
\end{cor}
\begin{pf} Again, all the arguments needed for  this proof are already in \cite{Bcrit}. 
That is, we prove exactly as in \cite[Prop.\ 4.3 and Cor.\ 4.4]{Bcrit} that the module
$\Symb_{\Gamma,C}^{\pm}(\D_{k}^\dag(L))_{(x)}$ is free of some rank $r$
over the algebra $\TT_{x,k}$ of the connected component of $x$ in the fiber of the weight map $w$ at $k$. The same argument shows that 
$S_{k+2}(\Gamma,L)_{(x)}$ is free of some rank $r'$ over $\TT_{x,k}$ and Proposition~\ref{compeigencurveBCM}  warrants that $r=r'$. The arguments of \cite[Prop. 4.5]{Bcrit}
then show that $r = \dim {S_{k+2}(\Gamma,L)[y]}$ for any $y$ satisfying our hypotheses (in {\it loc.\ cit}, we argue further that this dimension is $1$ under a newness hypothesis that we are not assuming here).
Finally, we note that $\dim \Symb_{\Gamma,C}^{\pm}(\D_{k}^\dag(L))[x]$ has dimension $r$ by the arguments of \cite[Theorem 4.7 and Cor. 4.8]{Bcrit}.
\end{pf} 

In order to apply the above corollary, we will need to know that $\CC_{\Gamma,\BCM}$ is smooth at certain points.
This is guaranteed by the following result.

\begin{theorem} \label{t:smootheigen}
Let $x$ be a point of $\CC_{\Gamma,\BCM}$ corresponding to the critical-slope refinenement of a decent new 
Eisenstein series $f$ of level $M$. We assume that $N/M$ is relatively prime to $M$ and that none of its prime divisors are bad for $f$ (cf.\ Def.~\ref{c:bad}). Then $\CC_{\Gamma,\BCM}$ is smooth at $x$.
\end{theorem}
\begin{pf} The case $N=M=1$ is the main theorem of \cite{BCsmooth}. 
The cohomological method used there
extends to the general case under the  hypothesis of the theorem. Let us comment briefly on the important
assumption that no prime divisors of $N/M$ are bad for $f$.

Suppose that $f = E_{k+2,\psi,\tau}$. Consider the subspace $H \subset H^1(G_\Q, \psi^{-1}\tau \omega^{k+1})$
 parameterizing extensions of 1 by $\psi^{-1} \tau \omega^{k+1}$ that are unramified outside $Np$
 and crystalline at $p$. The duality formula for Selmer groups implies that the dimension of $H$ is $1$ 
 plus the number of bad primes dividing $N/M$.
 When $H$ has dimension 1, the proof of the smoothness of the eigencurve goes exactly as in \cite{BCsmooth}.
Full details are given in \cite{Bcourse}.
\end{pf}

\section{Overconvergent boundary modular symbols} \label{s:overconvergent}

In this section, we denote by $C_0$ the set of cusps that are $\Gamma_0(p)$-equivalent to $0$. Thus, 
$C_0 = \Q \cap \Z_p$, the set of rational numbers that are $p$-integral. For  a cusp $c$, we denote by $\Gamma_c$ the stabilizer of $c$ in 
$\Gamma$.  Other notations used above (especially those of \S\ref{s:notation}) remain in force.

The monoid we called $S'$ in \S\ref{s:notation} is contained in $S_0(p)$. Therefore, as explained {\it loc.\ cit.},  the modules of full modular symbols $\Symb_{\Gamma}(\DD^\dag_\kappa[r](L))$ (or the same without the $\dag$, or the same with boundary modular symbols) have a natural action of the Hecke algebra $\HH$, the involution $\iota$, and in addition actions of 
 $U_\ell$ for $\ell \mid N$. We also have maps 
 \[ V_t: \Symb_{\Gamma_1(M) \cap \Gamma_0(p)}(\DD^\dag_\kappa[r](L))  \rightarrow \Symb_{\Gamma}(\DD^\dag_\kappa[r](L)) \] for integers $M,t$ such that $Mt \mid N$.

\subsection{Properties of  overconvergent boundary modular symbols}

We fix a finite extension $L$ of $\Q_p$, and for simplicity of notation, we shall write 
$\D_k$  for $\D_k(L)$, $\D_k^\dag$ for $\D_k^\dag(L)$, etc. 
 
\begin{lemma} \label{bsymb1} Let $k \in \Z$. If $\Phi \in \BSymb_\Gamma(\D_k^\dag[1])$ and $c$ is a cusp not in $C_0$, 
then $\Phi(c)=0$.
\end{lemma}
\begin{pf} 
We  need to show that $(\D_k[1]^\dag)^{\Gamma_c}=0$ whenever $c$ is a cusp in the $\Gamma_0(p)$ class of $\infty$. 
To see this, let $c = \gamma \infty$ with $\gamma \in \Gamma_0(p)$, and $\mu \in (\D_k^\dag[1])^{\Gamma_c}$.  Then $\mu_{|\gamma^{-1}}$ is invariant by $\gamma \Gamma_c \gamma^{-1} \supset \Gamma_\infty^n$ for some positive $n$ and we are reduced to the case $c=\infty$. 
In this case, $\mat{1 & n \\ 0 & 1 } \in \Gamma_\infty^n$, so $\mu(f(z+n))=\mu(f)$
for all functions $f$. Such a distribution is $0$ by \cite[Prop 3.1]{stevenspollack2}.
\end{pf}

\begin{prop} \label{bsymb2}
Let $j \geq 0$ be an integer. Let $\Phi \in \BSymb_\Gamma(\D_k^\dag)$ be a generalized eigenvector for $U_p$ corresponding to an eigenvalue $\lambda \neq 0$. We assume
that if $\lambda$ is of the form $\epsilon p^i$ where $\epsilon$ is a root of unity and $i \geq 0$ is an integer, then $i=j$. 
 Then for $c \in C_0$,
$\Phi(c)$ is a multiple of the distribution $\delta_{c,j}: f \mapsto \frac{df}{dz^{j}}(c)$.
\end{prop}
\begin{pf}
Since $\lambda \neq 0$, we have $\Phi \in \BSymb_\Gamma(\D_k^\dag[1])$
by Prop.~\ref{symbspace}(i). In particular, the preceding lemma applies and we have $\Phi(c')=0$ for any $c' \not \in C_0$.

Now let $c \in C_0$.
By hypothesis, we have $\Phi_{|(U_p - \lambda)^r} = 0$ for some integer $r$. The case $r=0$ is trivial, and we argue by induction over $r$, assuming the claim true for $r-1$.
For any integer $n \geq 1$, the operator $(U_p^n - \lambda^n) (U_p-\lambda)^{r-1}$ sends $\Phi$ to $0$.  Hence by the induction hypothesis,  we have 
\begin{eqnarray} \label{Upn} 
 ( \Phi_{|U_p^n} - \lambda^n \Phi)(c) = x \delta_{c,j}
 \end{eqnarray}
for some scalar $x \in L$.
 
Let us write $c=\sum_{i=0}^\infty a_i p^i$, with $0 \leq a_i \leq p-1$. Since $c$ is rational,
the sequence $(a_i)$ is eventually periodic: there exist two integers $m \geq 0$ and $n>0$ such that for all $i \geq m$, $a_{i+n}=a_i$. Let us assume for now that $m=0$. From (\ref{Upn}) and the definition of $U_p^n$, one gets
$$ \sum_{a=0}^{p^n-1} \Phi( (c-a)/p^n) )_{|\mat{1 & -a \\ 0 & p^n}} = x \delta_{c,j} + \lambda^n \Phi(c).$$
In the sum on the LHS, only one term is not 0, namely the term where $a=\sum_{i=0}^{n-1} a_i p^i$, because otherwise 
$(c-a)/p^n$ is not $p$-integral.
For this $a$, we have $c = (c-a)/p^n$, so we get
$$ \Phi(c)_{|\mat{1 & -a \\ 0 & p^n}} = x \delta_{c,j} + \lambda^n \Phi(c).$$
This implies that for  $f \in \Ac^\dag[1]$,
$$ \Phi(c)(f(p^n z +a) - \lambda^n f(z)) = x  \frac{df}{dz^{j}}(c).$$
Looking at the Taylor expansion at $c$ and using the assumption that $\lambda^n$ is not of the form 
$p^{in}$ for $i \geq 0$ an integer, except perhaps for $i=j$, one easily sees that the map 
\begin{align*}
\ker \delta_{c,j} \rightarrow& \  \ker \delta_{c,j}, \\
 f \mapsto& \ f(p^n z + a) - \lambda^n f(z)
 \end{align*}
  is surjective. Therefore, $\Phi(c)$ is $0$ on $\ker \delta_{c,j}$ and hence is a multiple of $\delta_{c,j}$.
It remains to treat the case $m>0$, but this can be reduced to the case $m=0$ using
(\ref{Upn}) for $m+n$ instead of $n$. Hence we have completed the induction step and the proof of the proposition.
\end{pf}
\begin{remark} Expressing the fact that $\Phi(c) \in (\D_k^\dag)^{\Gamma_c}$  for $c \in C_0$
 easily leads to the conclusion that if in the proposition above, $j$ is not $0$ or $k+1$, then 
 $\Phi(c)$ is actually $0$. We shall not use this result.
\end{remark}
\begin{definition} For $k \in \Z$ we call $\BSymb_{\Gamma}(\D_k^\dag)^\ord$ 
(resp.\ $\BSymb_{\Gamma}(\D_k^\dag)^\crit$ in the case $k \geq 0$) the linear span
of vectors in $\BSymb_{\Gamma}(\D_k^\dag)$ that are killed by $P(U_p)$ for some non-zero polynomial $P \in L[X]$
all of whose all roots in $\bar L$ have $p$-adic valuation $0$ (resp. $k+1$).
\end{definition}
\begin{theorem} \label{thetakisom} For $k \geq 0$, the map $$\Theta_k\colon  \BSymb_{\Gamma}(\D_{-2-k}^\dag)^\ord(k+1)
\rightarrow \BSymb_{\Gamma}(\D_k^\dag)^\crit$$ is a Hecke-equivariant isomorphism.
\end{theorem}
\begin{pf} The injectivity of $\Theta_k\colon \BSymb_{\Gamma}(\D_{-2-k}^\dag)(k+1)
\rightarrow \BSymb_{\Gamma}(\D_k^\dag)$ follows from (\ref{fundexsymb}) since $\BSymb_{\Gamma}$ is 
 a left-exact functor. It is clear that $\Theta_k$ maps 
$\BSymb_{\Gamma}(\D_{-2-k}^\dag)^\ord$ into $\BSymb_{\Gamma}(\D_{k}^\dag)^\crit$ since $U_p$-eigenvalues are multiplied by $p^{k+1}$. These spaces are finite-dimensional by Prop.~\ref{symbspace}(i), hence after extending $L$ if necessary we can assume that they are generated by generalized eigenvectors for $U_p$ with eigenvalue $\lambda$ of $p$-adic order $0$ and $k+1$, respectively. Hence by Lemma~\ref{bsymb1} and Proposition~\ref{bsymb2}, any vector $\Phi$ in the target has the form 
\[ \Phi(c) = \begin{cases} 0 & \text{ if  } c \not \in C_0 \\
\tau(c) \delta_{c,k+1} & \text{ if } c \in C_0, \end{cases} \]
for some $\tau(c) \in L$.
Define $\Phi' \in \Hom(\Delta,\D_{-2-k})$ by 
\[ \Phi'(c) = \begin{cases} 0 & \text{ if  } c \not \in C_0 \\
\tau(c) \delta_{c,0} & \text{ if } c \in C_0. \end{cases} \]
 Then clearly $\Theta_k \Phi' = \Phi$. Since $\Theta_k$ is $\Gamma$-equivariant and injective, $\Phi'$ is
$\Gamma$-invariant, i.e.\ is an element of $ \BSymb_{\Gamma}(\D_{-2-k}^\dag)^\ord$.
This proves the surjectivity of $\Theta_k$.
\end{pf}

\begin{remark} \label{r:signboundary} We observe from the proof of the theorem that any critical boundary modular symbol of trivial nebentypus and square free level has sign $-1$. Indeed, such a modular symbol $\Phi$ has $\Phi(c)$ a multiple of $\delta_{c,k+1}$ with $k$ even, and $c$ and $-c$ are in the same class modulo $\Gamma_0(M)$ for any square free $M$.
\end{remark}

The following is a version of Stevens' control theorem for boundary symbols.
\begin{prop} \label{controlboundary} Assume $k \geq 0$. The map $\rho_k\colon \BSymb_\Gamma(\D_k^\dag)^\ord \rightarrow \BSymb_\Gamma(\V_k)^\ord$ is an isomorphism.
\end{prop}
\begin{pf} 
The map  $\rho_k\colon \BSymb_\Gamma(\D_k^\dag)^\ord \rightarrow \BSymb_\Gamma(\V_k)^\ord$ is injective since it is the restriction of  $\rho_k\colon \Symb_\Gamma(\D_k^\dag)^\ord \rightarrow \Symb_\Gamma(\V_k)^\ord$, which is injective by Stevens' control theorem (Proposition~\ref{symbspace}(iii)). For the surjectivity, we just have to check that if $\phi \in \BSymb_\Gamma(\V_k)^\ord$, then for each cusp $c$, $\phi(c) \in (\V_k)^{\Gamma_c}$ can be lifted to a distribution in $(\D_k^\dag)^{\Gamma_c}$. From Corollary~\ref{c:ordphi} and the succeeding comments, we see that
$\phi(c) \in \V_k$ is $0$ when $c \not \in C_0$, and $\phi(c)$ is a multiple of the linear form $P \mapsto P(c)$ when $c \in C_0 \subset \Z_p$. Such a linear form can be lifted to the  invariant distribution $f \mapsto f(c)$. 
\end{pf}
 
\subsection{A basis of the space of ordinary overconvergent boundary modular symbols}

In this subsection, $L$ is a finite extension of $\Q_p$ that contains the $\varphi(N)$th roots of unity.

For $k \in \Z$, $M$ a positive integer that divides $N$, $u,v$ two relatively prime integers such that
$p \nmid v$, define 
$\Phi_{k,u,v} \in \BSymb_{\Gamma_1(M) \cap \Gamma_0(p)}(\D^\dag_k(L))$  as follows.
$\Phi_{k, u, v}$ is supported on the $\Gamma_1(M) \cap \Gamma_0(p)$-orbit of $u/v$, and is defined on 
that orbit
by the same formula as (\ref{e:phiuvdef}), with $P(z) \in \P_k$ replaced by $f(z) \in \Ac^\dag_k$, i.e.\ 
\begin{equation} \label{e:Phiuvdef}
 \Phi_{k, u,v}\left(\gamma\left(\frac{u}{v}\right)\right)(f(z)) = f\left(\gamma\left(\frac{u}{v}\right)\right) \cdot (cu + dv)^k, 
 \end{equation}
 where
\[  \gamma = \mat{a & b \\  c  & d} \in \Gamma_1(M) \cap \Gamma_0(p). \]

For any decomposition $M=QR$ and any two Dirichlet characters 
$\psi, \tau$ with conductors $Q, R$, respectively, such that $\psi\tau(-1)=(-1)^k$, define $\Phi_{k,\psi,\tau}$ by the same formula (\ref{e:phipsitaus}) as $\phi_{k,\psi,\tau}^p$ with $\phi_{k,x,Qy}$ replaced by $\Phi_{k,x,Qy}$, i.e. 
\begin{equation} \label{e:Phipsitau} \Phi_{k, \psi, \tau} = 
\sum_{\genfrac{}{}{0pt}{}{x \!\!\!\! \pmod{Q}}{(x,Q)=1}}
\sum_{\genfrac{}{}{0pt}{}{y \!\!\!\! \pmod{R}}{(y,Rp)=1}}
\psi^{-1}(x) \tau(y) \Phi_{k, x, Qy}.\end{equation}
Then $\Phi_{k,\psi,\tau} \in \BSymb_{\Gamma_1(M) \cap \Gamma_0(p)}(\D_k^\dag(L))$, and
we clearly have $\rho_k(\Phi_{k, \psi, \tau}) = \phi_{k, \psi, \tau}^{p}$ for $k \ge 0$.

The symbol $\Phi_{k,\psi,\tau}$ is a Hecke eigenvector with the same eigenvalues as 
$\phi_{k,\psi,\tau}^{p}$, namely:
\begin{itemize}
\item[(i)]
for $T_\ell \  (\ell \nmid Np)$ (resp. $U_\ell$ for $\ell \mid N$) with  eigenvalues $\psi(\ell) + \tau(\ell)\ell^{k+1}$,
\item[(ii)] for $\langle d \rangle \ (d \in (\Z/n\Z)^\ast)$  with
eigenvalues $\psi(d)\tau(d)$,
\item[(iii)] for $U_p$ with eigenvalue $\psi(p)$,
\item[(iv)] and for $\iota$ with eigenvalue $\psi(-1)$.
\end{itemize}
In particular, $\Phi_{k,\psi,\tau} \in \BSymb_{\Gamma_1(M) \cap \Gamma_0(p)}(\D_k^\dag(L))^\ord$.

As usual we set
$\Phi_{k,\psi,\tau,t}=(\Phi_{k,\psi,\tau})_{|V_t}$.
For $t$ such that $QRt \mid N$, we have 
$\Phi_{k,\psi,\tau,t} \in \BSymb_\Gamma(\D_k^\dag(L))^\ord$.
Moreover, \begin{equation} \label{e:rhokphi}
\rho_k (\Phi_{k,\psi,\tau,t}) = \phi_{k,\psi,\tau,t}^p 
\end{equation} for $k \geq 0$.

\begin{theorem} \label{basisboundover} Let  $k \in \Z$, $k \neq -1$.
Let $\psi, \tau$ be Dirichlet characters of conductors $Q, R$, respectively, 
such that $\tau\psi(-1) = (-1)^k$ and  $M= QR$ divides $N$.  Suppose that
$N/M$ is squarefree, relatively prime to $M$, and such that no primes dividing $N/M$
are bad for $E_{k+2, \psi, \tau}$.
Then the symbols $\Phi_{k,\psi,\tau,t}$ for positive integers $t$ dividing $N/M$
form a basis of the space 
$ \BSymb_{\Gamma}(\D_{k}^\dag(L))^\ord[E_{k+2, \psi, \tau}]$.
\end{theorem}
\begin{pf} For $k \geq 0$, this follows from Proposition~\ref{controlboundary}, Corollary~\ref{c:ordphi}, and (\ref{e:rhokphi}).
We will conclude the result for negative integers by interpolation.  To this end, 
fix  $k < 0$, and let $X=\sp A$ be  a closed ball with radius in $p^{\Q}$ in the
connected component of $\W$ containing $k$. Then $A$ is a principal affinoid domain.
The module $\BSymb_\Gamma(\D_K^\dag(A))^\ord$ is free of finite rank by Prop.~\ref{symbmodule}(ii).  It follows that the fibers at all points $\kappa \in X(\Q_p)$ of that module (namely $\BSymb_\Gamma(\D_K^\dag(A))^\ord \otimes_{A,\kappa} \Q_p$ where the implicit map $A \rightarrow \Q_p$ is the one defined by $\kappa$) have the same dimension. 

We claim that when $\kappa$ is an integer $k' \in \Z$, the natural map
 $$\BSymb_\Gamma(\D_K^\dag(A))^\ord \otimes_{A,k'} \Q_p \rightarrow \BSymb_\Gamma(\D_{k'}^\dag(\Q_p))^\ord$$
 is a Hecke-equivariant isomorphism.  The injectivity of this map follows formally form the left exactness
of the functor $\BSymb_\Gamma$. For the surjectivity, we need to know that if
$\Phi \in \BSymb_\Gamma(\D_{k'}^\dag(\Q_p))^\ord$, then for all cusps $c$,
the distribution $\Phi(c)$ may be extended to a $\Gamma_c$-invariant distribution in $\D_K^\dag(A)$.
 But $\Phi(c)$ is $0$ if $c \not \in C_0$ by Lemma~\ref{bsymb1}, and $\Phi(c)$ is a multiple of the 
 ``evaluation at $c$" distribution  $\delta_{0,c}$  if $c \in C_0 \subset \Z_p$ by Prop.~\ref{bsymb2}. In both cases, the distributions can obviously be lifted.
This proves the claim.

Choosing a $k' > 0$ in $X$, we get that 
 $$\dim \BSymb_\Gamma(\D_{k'}^\dag(\Q_p))^\ord[E_{k+2, \psi, \tau}] = \dim \BSymb_\Gamma(\D_{k}^\dag(\Q_p))^\ord[E_{k+2, \psi, \tau}].$$ 
  This dimension is the cardinality of the family of elements $\Phi_{k,\psi,\tau,t}$
 in the statement of the theorem, namely, the number of divisors $t$ of $N/M$.  With $\psi$,  $\tau$, and $k \neq -1$  fixed, the symbols $\Phi_{k,\psi,\tau,t}$ for $t \mid N/M$
 are linearly
 independent; this is easily seen as in Corollary~\ref{c:ordphi}, since the vectors
 $(\Phi_{k,\psi,\tau})_{|O_\ell}$ and  $(\Phi_{k,\psi,\tau})_{|C_\ell}$ for primes $\ell \mid N/M$ have distinct $U_\ell$-eigenvalues.
  Therefore, the symbols $\Phi_{k,\psi,\tau,t}$ form a basis of $\BSymb_\Gamma(\D_{k}^\dag(L))^\ord$ as desired.
 \end{pf}

\subsection{Proof of the $-\epsilon(f)$-part of Theorem~\ref{t:main}}

\begin{prop} \label{p:minuspepsart} Let $f$ be any new Eisenstein series, $f_\beta$ its critical refinement, and $\mu_{f_\beta}^{-\epsilon(f)}$
defined as in the introduction. Then $\mu_{f_\beta}^{- \epsilon(f)} = 0$
unless $f$ is a normal Eisenstein series $E_{k+2, \psi, \tau}$ with $\tau = 1$
or $f$ is  an exceptional  Eisenstein series  $E_{2,\ell}$. 
In the latter cases, the distribution $\mu_{f_\beta}^{-\epsilon(f)}$ has support $\{0\}$.
 In all cases, we have
$$L_p(f_\beta,\sigma)=0 \text{\ \ \ if $\sigma(-1)=-\epsilon(f)$}.$$
\end{prop}

\begin{pf}
The assertion about the $p$-adic $L$-function follows from the assertions about $\mu_{f_\beta}^{-\epsilon(f)}$,
since the Mellin transform of a distribution with support contained in $p \Z_p$ is $0$. 

We only treat the normal case, leaving the simple modifications for the exceptional case to the reader.  Let 
$f=E_{k+2,\psi,\tau}$, with $\psi$ and $\tau$ primitive of conductor $Q$ and $R$, respectively.
Then $f$ is a newform of level $N := QR$, and $f_\beta$ is a form of level $\Gamma=\Gamma_1(N) \cap \Gamma_0(p)$.
 We have $\epsilon(f)=\psi(-1)$ (cf.\ Prop.~\ref{p:boundsymb}). By definition, $\Phi_{f_\beta}^{-\epsilon(f)}$ is a generator of the one-dimensional
space $\Symb_{\Gamma}^{- \epsilon(f)}(\D_k^\dag)[f_\beta]$. We observe that the modular symbol 
$\Theta_k \Phi_{-2-k,\tau,\psi}$ (note the inversion in the order of $\psi$ and $\tau$) has the same eigenvalues 
for the Hecke operators $T_\ell$, $U_p$ and $\langle d \rangle$ as $f_\beta$, 
and has sign $\tau(-1) (-1)^{k+1} = - \psi(-1) = - \epsilon(f)$.
Therefore, up to a non-zero scalar,
$$\Phi_{f_\beta}^{-\epsilon(f)} = \Theta_k \Phi_{-2-k,\tau,\psi}$$
and
$\mu_{f_\beta}^{-\epsilon(f)}$ is the $(k+1)$-derivative of the measure  $\Phi_{-2-k,\tau,\psi}(\{\infty\} - \{0\})$.
This measure is readily seen by its definition (\ref{e:Phipsitau}) to be $0$ unless  $\tau=1$, in which case it is a multiple of the ``evaluation at $0$" measure. This completes the proof.
\end{pf}

\begin{remark} It follows from the proof above that $\Phi_{f_\beta}^{- \epsilon(f)}$ is a boundary modular symbol.
It is easy to see (cf. Remark~\ref{r:signboundary} above for certain cases or Lemma~\ref{lemmabsymb0} below in general) that on the contrary $\Phi_{f_\beta}^{\epsilon(f)}$ is not a boundary modular symbol.
\end{remark}

\section{Families of $p$-adic $L$-functions of Eisenstein series: proof of the main theorem in the normal case} \label{s:families}

In this section we prove the main theorem in the normal case.
As in the introduction,  we fix $\tau$ a Dirichlet character of conductor $R$, $\psi$ a Dirichlet character of conductor $Q$, and a prime number $p$ not dividing $QR$.
We will also choose a multiple $N$ of $QR$ not divisible by $p$;  we will have either $N=QR$, $N=QR \ell$,
or $N=QR\ell_1\ell_2$ for primes $\ell, \ell_1, \ell_2$ not dividing $QRp$.
We use the notations defined in \S\ref{s:notation} for our given choice of $N$.

\subsection{The $p$-adic $L$-functions of an ordinary Eisenstein series}
\label{s:61}

Consider the Eisenstein series $E_{k+2,\tau,\psi,t}$ of level $\Gamma_1(N)$
where $t$ is a positive integer such that $QRt \mid N$.  (As noted earlier, the order of $\tau$ and $\psi$ has been switched from the introduction for our intended application; namely, the  operator $\Theta_k$ sends $E_{-k, \tau, \psi}^\ord$ to $E_{k+2, \psi, \tau}^\crit$.)
Like any eigenform of that level, $E_{k+2,\tau,\psi,t}$ has two $p$-refinements, which are eigenforms for $\HH$ of level $\Gamma = \Gamma_1(N) \cap \Gamma_0(p)$. They are:
\begin{eqnarray*} E_{k+2,\tau,\psi,t}^\ord &=& E_{k+2,\tau,\psi,t} - \psi(p) p^{k+1} E_{k+2,\tau,\psi,pt}, \\  E_{k+2,\tau,\psi,t}^\crit &=& E_{k+2,\tau,\psi,t} - \tau(p) E_{k+2,\tau,\psi,pt}.
\end{eqnarray*}
Both forms have the same eigenvalues as $E_{k+2,\tau,\psi,t}$ for the Hecke operators $T_\ell$ ($\ell \nmid Np$) and the Diamond operators, but they have the following eigenvalues for $U_p$:
\begin{eqnarray*} U_p E_{k+2,\tau,\psi,t}^\ord &=& \tau(p) E_{k+2,\tau,\psi,t}^\ord  \\
U_p  E_{k+2,\tau,\psi,t}^\crit &=& \psi(p) p^{k+1} E_{k+2,\tau,\psi,t}^\crit.
\end{eqnarray*}
Since the $U_p$-eigenvalue of $E_{k+2,\tau,\psi,t}^\ord$ is a $p$-adic unit, this form is ordinary.

We consider a linear combination 
\begin{equation} \label{e:fdef}
f = \sum_{t, QRt | N} c_t  E_{k+2,\tau,\psi,t}^\ord
\end{equation}
where the $c_t$ are algebraic numbers. The form $f$ is an eigenform of level $\Gamma$. Let us assume that $f$ is $C$-cuspidal, so we can consider the partial
modular symbol $\phi_f \in 
\Symb_{\Gamma,C}(\V_k(K))$ defined in \S\ref{s:ccuspidal}.
Here $K \subset \overline \Q$ is the number field of Proposition~\ref{symbeisgood}.
We consider
the part of $\phi_f$ with sign $-\tau(-1)$, and using the embedding $\overline{\Q} \subset \overline{\Q}_p$ we
view the coefficients as living in $\V_k(L)$, where $L$ is the finite extension of $\Q_p$ generated by $K$: 
\[ \phi_f^{-\tau(-1)} \in \Symb_{\Gamma, C}(\V_k(L)). \] 
By Prop~\ref{symbspace}(iii), there is a unique $\Phi_f \in \Symb^{-\tau(-1)}_{\Gamma,C}(\D_k(L))$ such that \[ \rho_k(\Phi_f)=\phi_f^{-\tau(-1)}. \] 

\begin{prop} \label{propLpf1} Let $f$ be as in $(\ref{e:fdef})$ such that $f$ is $C$-cuspidal.
If $\sigma \in \WW(\C_p)$ and $\sigma(-1) = - \tau(-1)$, then
\begin{equation} \label{decLp} L_p(\Phi_f,\sigma) = \frac{G(\tau)}{2R} \sum_{t, \ QRt \mid N} c_t  t^{-1} \sigma(Rt)^{-1} L_p(\psi,\sigma z^{-k}) L_p(\tau,\sigma z).
\end{equation}
\end{prop}
\begin{pf} 
We first recall the interpolation properties of the Kubota--Leopoldt $p$-adic $L$-functions
$L_p(\nu,\sigma)$ for $\nu$ a Dirichlet character of conductor prime to $p$.
If $\chi$ is a finite order character of $\Z_p^\ast$ of  conductor $p^n$, and $m$  is an integer, we have: 
\begin{equation} \label{dirichletintera}
L_p( \nu, \chi z^m ) = L(\nu \chi^{-1}, m)  
\end{equation}
if $\nu\chi(-1) = (-1)^{m+1}$ and $m \le 0$.  Also, 
\begin{align} \label{dirichletinterb} 
L_p( \nu, \chi z^m ) &= \frac{2\Gamma(m) (\cond(\nu)p^n)^m \nu\chi^{-1}(-1)}{(2 \pi i)^m G(\nu \chi^{-1})} L(\nu \chi^{-1},m) \\
&= L(\nu^{-1}\chi, 1-m) \label{dirichletinterc}
\end{align}
if $\nu\chi(-1) = (-1)^m$ and $m \ge 1$.  In (\ref{dirichletinterb}), $G(\nu\chi^{-1})$ denotes a Gauss sum.
Note that if $\chi = 1$,  these three equations must be altered by removing the Euler factors at $p$ in 
(\ref{dirichletintera}) and (\ref{dirichletinterc}); in other words, the right side of (\ref{dirichletintera}) should be 
multiplied by $(1 - \nu(p)p^{-m})$, and the right sides of (\ref{dirichletinterb}) and (\ref{dirichletinterc}) should be
multiplied by $(1 - \nu^{-1}(p)p^{m-1})$.

Let us assume that $\sigma$ is a special character of the form $z \mapsto  \chi(z) z^j$ where $\chi$ is a finite order character of $\Z_p^\ast$ of conductor $p^n$ and $0 \leq j \leq k$,
and that $\sigma(-1) = -\tau(-1)$.
A simple and standard computation using the fact that $\rho_k(\Phi_f)=\phi_f^{-\tau(-1)}$, $U_p \Phi_f = \tau(p) \Phi_f$, and that $\phi_f$ satisfies (\ref{formulaLf2})
allows one to compute  $L_p(\Phi_f,\chi z^j)$. 
If $\chi$ is non-trivial, one arrives at the interpolation formula (\ref{noncritinter}) for $f$:
 $$L_p(\Phi_f,\chi z^j)= \frac{p^{n(j+1)} j!}{\tau(p)^n (-2 \pi i)^{j+1} G(\chi^{-1})} L(f,\chi^{-1},j+1),$$
where $G(\chi^{-1})$ is the Gauss sum. 
Using Lemma~\ref{Leisensteinlemma}, we obtain 
\begin{equation} \label{Lp1}  L_p(\Phi_f,\chi z^j)= \!\!\!\sum_{t,\ QRt|N} \!\!\!\!\!
c_t \frac{p^{n(j+1)} j!  t^{-j-1} \chi^{-1}(t)}{\tau(p)^n (-2 \pi i)^{j+1} G(\chi^{-1})} L(\tau \chi^{-1},j+1) L(\psi \chi^{-1},j-k).\end{equation}

To express each term of this sum in terms of $p$-adic $L$-functions, we apply (\ref{dirichletintera}) to $\nu = \psi$ and $m = j-k$. This is possible since $m \leq 0$ and 
\[  \psi\chi(-1)(-1)^{j-k + 1}
 = \psi \sigma(-1) (-1)^{-k+1} = - \psi \tau(-1) (-1)^{-k+1} = 1, \]  and we obtain
 \begin{eqnarray} \label{dir3} 
 L(\psi \chi^{-1},j-k) = L_p( \psi, \chi z^{j-k} ).
\end{eqnarray}
Next we apply (\ref{dirichletinterb}) to $\nu = \tau$ and $m = j+1$. This is possible since $m \geq 1$ and 
\[ \nu\chi(-1)(-1)^m = \tau\chi(-1)(-1)^{j+1} = (-\tau(-1))^2 = 1,\]
and we obtain
 \begin{eqnarray} \label{dir4}
 L(\tau \chi^{-1},j+1) = \frac{(2 \pi i)^{j+1} G(\tau\chi^{-1}) \tau\chi(-1)} {2 j! (Rp^n)^{j+1}} L_p( \tau, \chi z^{j+1}).
\end{eqnarray}
Plugging in (\ref{dir3}) and (\ref{dir4}) and using $G(\tau \chi^{-1}) = \tau(p^n)\chi^{-1}(R) G(\tau) G(\chi^{-1})$,
we find
$$  L_p(\Phi_f,\chi z^j)=  \frac{G(\tau)}{2}\sum_{t} c_t \chi(Rt)^{-1} (Rt)^{-j-1}  L_p( \psi, \chi z^{j-k} ) L_p( \tau, \chi z^{j+1}),$$ 
which is (\ref{decLp}) for the characters $\sigma =\chi z^j$ satisfying $ 0 \leq j \leq k$
and $\sigma(-1) = -\tau(-1)$. There are infinitely many such characters on each
component of $\WW$ of sign $-\tau(-1)$, and  both sides of (\ref{decLp}) are bounded functions. The result follows since a non-zero bounded analytic function on an open ball has at most finitely many zeros.
\end{pf}

Let us now choose specific constants $c_t$ such that the resulting form $f$ is $C$-cuspidal.  The form $f$ is cuspidal at $\infty$ automatically if $\tau \neq 1$, and if $\sum_{t} c_t =0$ when $\tau=1$.  Similarly, the form $f$ is cuspidal at $0$ automatically if $\psi \neq 1$, and if $\sum c_t t^{-k-2} = 0$ when $\psi=1$.  It will be convenient to choose $N$ minimal necessary to achieve $C$-cuspidality.  We therefore propose:

\begin{definition} \label{d:Ndef} Define $N$ and a modular form $F_{k+2}$ of level $\Gamma$ as follows.
\begin{itemize} \item If $\tau \neq 1$, $\psi \neq 1$, let $N = QR$ and $F_{k+2} = E_{k+2,\tau,\psi}^{\ord}.$
\item If $\tau = 1$  but $\psi \neq 1$, choose a prime $\ell_1 \nmid QRp$ such that $\psi(\ell_1) \neq 1$, let $N = QR\ell_1$,
and let $F_{k+2} = E_{k+2,\tau,\psi}^{\ord} - E_{k+2,\tau,\psi, \ell_1}^{\ord}.$
\item If $\psi = 1$  but $\tau \neq 1$, choose a prime $\ell_2 \nmid QRp$ such that $\tau(\ell_2) \neq 1$, let $N = QR\ell_2$,
and let $F_{k+2} = E_{k+2,\tau,\psi}^{\ord} - \ell_2^{k+2} E_{k+2,\tau,\psi, \ell_2}^{\ord}.$
\item If $\tau = \psi = 1$ in the normal case (so $k > 0$), choose primes $\ell_1$ and $\ell_2$
not dividing $QRp$, let $N = QR\ell_1\ell_2$, and let
\[ 
F_{k+2} = E_{k+2,\tau,\psi}^\ord - E_{k+2,\tau,\psi,\ell_1}^\ord- \ell_2^{k+2} E_{k+2,\tau,\psi,\ell_2}^\ord + \ell_2^{k+2} E_{k+2,\tau,\psi,\ell_1\ell_2}^\ord. \]
\end{itemize}
\end{definition}

Note that by construction, in the normal case the assumptions of Corollary~\ref{c:ordphi} and Theorem~\ref{basisboundover} hold for our choice of $N$, i.e.\ no primes dividing $N/(QR)$ are bad for $E_{k+2, \psi, \tau}$.

From Proposition~\ref{propLpf1} and the succeeding comments, we obtain:

\begin{prop} \label{p:lfk}
The Eisenstein series $F_{k+2}$ is $C$-cuspidal,
and for $\sigma \in \W(\C_p)$ with $\sigma(-1) = -\tau(-1)$ we have  $$L_p(F_{k+2},\sigma)  = m(k,\sigma)
 L_p(\psi,\sigma z^{-k}) L_p(\tau,\sigma z),$$
 where 
 \begin{equation} \label{e:mdef}
  m(k, \sigma) = \frac{G(\tau)}{2R \sigma(R)} (1- \sigma(\ell_1)^{-1} \ell_1^{-1}) (1- \sigma(\ell_2)^{-1} \ell_2^{k+1}),
  \end{equation}
with the understanding that the factor associated to $\ell_1$ or $\ell_2$ does not occur in $(\ref{e:mdef})$ if the corresponding
prime does not occur in the relevant case of Definition~\ref{d:Ndef}.
\end{prop}

\subsection{A family of partial modular symbols}
\label{s:62}

 In joint work with H.~Darmon, the second-named author proved that $\Phi_{F_{k+2}}$ belongs to a family of partial modular symbols over the weight space \cite{DD}.
If  $X=\sp A$ is an open affinoid subset of $\WW$,  $K \in \WW(A)$ is the canonical character,  $\Phi \in \Symb_{\Gamma,C}( \DD_K(A))$, and $k \in \Z$, then we write 
 $\Phi_k$ for the image of $\Phi \otimes 1$ in  
 $ \Symb_{\Gamma,C}( \DD_\kappa )(\Q_p) $ by the specialization map (\ref{specmap}).

\begin{theorem}[Dasgupta-Darmon] \label{t:DD}There exists a unique symbol \[ \Phi \in \Symb_{\Gamma,C}( \DD_K (A))\] such that for every integer $k \ge 0$, we have $\Phi_k = \Phi^{-\tau(-1)}_{F_{k+2}}$.
\end{theorem}
\begin{pf} 
 Let $M(\Z_p \times \Z_p^\times)$ denote the $\Gamma$-module of $\Z$-valued measures on $\Z_p \times \Z_p^\times$.
In \cite[Theorem 4.2]{DD}, it is proved that there exists a unique partial modular symbol $\mu \in \Symb_{\Gamma, C}(M(\Z_p \times \Z_p^\times))$ such that
for each homogenous polynomial $h(x,y) \in \Q[x,y]$ of degree $k$, we have
\[ \int_{\Z_p \times \Z_p^\times} h(x,y) d\mu(\{r\}- \{s\})(x,y) = \phi_{F_{k+2}}^{-\tau(-1)}(\{r\}- \{s\})(h(z, 1)). \]
It is then clear that
\[ \Phi(\{r\} - \{s\})(g(z)) := \int_{\Z_p \times \Z_p^\times} g(x/y)K(y) d\mu(\{r\}- \{s\})(x,y) \]
has the desired property.
\end{pf}

\begin{prop} \label{p:phiev} Let $\Phi$ be as in Theorem~\ref{t:DD}.
Assume that $\Z^{\ge 0} \cap X$ is Zariski-dense in $X$.
Then \begin{eqnarray*}
\Phi_{|\iota} &=& - \tau(-1) \Phi \\ 
\Phi_{|T_\ell} &=& ( \tau(\ell) + K(\ell) \ell \psi(\ell) )\Phi, \quad \ell \nmid Np \\
\Phi_{|U_p} &=& \tau(p) \Phi\\
\Phi_{|\diamond{a}} &=& \tau(a) \psi(a) \Phi, \quad a \in (\Z/N\Z)^\times.
\end{eqnarray*}
\end{prop}
\begin{pf} These results hold after applying the specialization map at any positive integer $k$.
Since $\Z^{\ge 0} \cap X$ is Zariski-dense, the proposition follows.
\end{pf}

\begin{theorem} \label{t:psidef} There exists $\Psi \in \Symb_{\Gamma,C}(\D_{-2-k}(L))$ such that
 \begin{eqnarray*}
 \Psi_{|\iota} &=& - \tau(-1) \Psi \\ 
\Psi_{|T_\ell} &=& ( \tau(\ell) + \ell^{-1-k} \psi(\ell) )\Psi, \quad \ell \nmid Np \\
\Psi_{|U_p} &=& \tau(p) \Psi\\
\Psi_{|\diamond{a}} &=& \tau(a) \psi(a) \Psi, \quad a \in (\Z/N\Z)^*
\end{eqnarray*}
and such that for $\sigma \in \W(\C_p)$ with $\sigma(-1) = -\tau(-1)$, we have $$L_p(\Psi,\sigma) = m(-2-k,\sigma)
 L_p(\psi,\sigma z^{k+2}) L_p(\tau,\sigma z).$$
\end{theorem}
\begin{pf} Choose for $X$ a closed ball around $-2-k$ in $\WW$. Then $\Z^{\ge 0} \cap X$ is  Zariski-dense in $X$.
Let $\Phi$ be as in Theorem~\ref{t:DD}.
Applying Proposition~\ref{p:phiev} to $X$, it follows that the specialization $\Psi:=\Phi_{-2-k}$ of $\Phi$ has the desired eigenvalues.

 Moreover, $L_p(\Phi)$ is an element of $\RR \hotimes A$, i.e.\ a two variable $p$-adic $L$-function $L_p(\Phi)(\kappa,\sigma)$ where $\kappa \in X(\C_p)$ and $\sigma \in \W(\C_p)$, such that $L_p(\Phi)(k,\sigma) = L_p(\Phi_k,\sigma)$ when $\kappa=k \in \Z$.
In particular, if $k \ge 0$, we have by Proposition~\ref{p:lfk}
\[ L_p(\Phi)(k,\sigma) =  m(k, \sigma) L_p(\psi,\sigma z^{-k}) L_p(\tau,\sigma z) \]
  whenever $\sigma(-1) = -\tau(-1).$
By interpolation, we have for all $\kappa$:
\[ L_p(\Phi)(\kappa,\sigma) = m(\kappa, \sigma)  L_p(\psi,\sigma \kappa^{-1}) L_p(\tau,\sigma z),\]
 where $m(\kappa, \sigma)$ has the obvious meaning (namely (\ref{e:mdef}) with $k$ replaced by $\kappa$).
 The result follows by specialization to $\kappa = -2 -k$.
\end{pf}

\begin{cor} \label{corphikcrit} There exists a modular symbol \[ \Phi_k^\crit \in 
\Symb_{\Gamma,C}^{\psi(-1)}(\D_k^\dag(L))[E_{k+2,\psi,\tau}] \] such that for $\sigma(-1) = \psi(-1)$ we have
$$ L_p(\Phi_k^\crit,\sigma) = \log^{[k+1]}(\sigma)n(k, \sigma)
 L_p(\psi,\sigma z) L_p(\tau,\sigma z^{-k}),$$ where
 \[ n(k, \sigma) = m(-2-k, \sigma z^{-k-1}). \]
\end{cor}
\begin{pf} Set $\Phi_k^\crit = \Theta_k\Psi$ where $\Psi$ is as in Theorem~\ref{t:psidef}, i.e.
\[ \Phi_k^\crit(\{r\} - \{s\})(g(z)) = \Psi(\{r\} - \{s\})\left( \frac{d^{k+1} g(z)}{dz^{k+1}} \right). \]
The map $\Theta_k$ multiplies the eigenvalues of the operators $\iota$, $T_\ell$, $U_p$, and $\langle a \rangle$ by
$(-1)^{k+1}, \ell^{k+1}, p^{k+1}$, and $1$.  Therefore, $\Phi_k^\crit$ is in the desired eigenspace of 
$\Symb_{\Gamma,C}(\D_k^\dagger(L))$: 
 \begin{eqnarray*}
  {\Phi_k^\crit}_{|\iota} &=& \psi(-1)  \Phi_k^\crit \\ 
 {\Phi_k^\crit}_{|T_\ell} &=& ( \psi(\ell) + \ell^{k+1} \tau(\ell) ) \Phi_k^\crit, \quad \ell \nmid Np \\
 {\Phi_k^\crit}_{|U_p} &=& \tau(p)p^{k+1}  \Phi_k^\crit\\
 {\Phi_k^\crit}_{|\diamond{a}} &=& \tau(a) \psi(a)  \Phi_k^\crit, \quad a \in (\Z/N\Z)^*.
\end{eqnarray*}

One easily calculates from the definition that in general one has
 \[ L_p(\Theta_k \Psi, \sigma) = \log^{[k+1]}(\sigma) L_p(\Psi, \sigma z^{-k-1}) \]
(see \cite[Lemma 2.9]{BCM}).   Note also that  \[ \sigma(-1) = \psi(-1) \Longrightarrow (\sigma z^{-k-1})(-1) = - \tau(-1).\]
The desired result now follows from the calculation of $L_p(\Psi, \sigma)$ in Theorem~\ref{t:psidef}.
\end{pf}
 
\subsection{End of the proof in the normal case}
\label{s:63}

Recall that  $\Gamma = \Gamma_1(N) \cap \Gamma_0(p)$.
Set \[ \Gamma' = \Gamma_1(QR) \cap \Gamma_0(p) \supset \Gamma. \]

Let $f =   {E_{k+2,\psi,\tau}}$  and $\beta = \tau(p) p^{k+1}$, and write $f_\beta$ for the
 critical refinement of $f$, also known as the ``evil" Eisenstein series.  The
 form $f_\beta$  is a modular form for $\Gamma'$.
 We view it as an old form for $\Gamma$.

Let $\Phi_{f_\beta}$ be a generator of the space $\Symb_{\Gamma'}^{\psi(-1)}(\D_k^\dag)[f_\beta]$ (which is one-dimensional by the main result of \cite{Bcrit}).  
 
 The following proposition is the analogue for modular symbols of  the linear independence of the modular forms $f_\beta(tz)$ for distinct integers $t$,  which is easily seen on $q$-expansions using the theory of newforms.

\begin{prop} \label{fourindependent} The vectors  $(\Phi_{f_\beta})_{|V_{t}}$ for $t \mid N/(QR)$ 
 are linearly independent in $\Symb_\Gamma^{\psi(-1)}(\D_k^\dag(L))[f_\beta]$.
\end{prop}
\begin{pf} 
We distinguish two cases, according to whether
the image of $\Phi_{f_\beta}$ by the map $\rho_k$ of (\ref{fundexseq}) is 0 or not.
In the latter case, the image $\rho_k(\Phi_{f_\beta})$ is (up to a nonzero scalar)
the boundary classical modular symbol $\phi_{k, \psi, \tau}^p \in \BSymb_{\Gamma'}(\V_k(L))$
by Corollary \ref{c:ordphi} (with $N$ replaced by $QR$).  Applying the corollary again (with $N$ playing its own role),
we see that the vectors $\rho_k((\Phi_{f_\beta})_{| V_t}) =  \phi_{k, \psi, \tau, t}^p$ are linearly independent,
and hence the vectors  $(\Phi_{f_\beta})_{|V_{t}}$ are as well.

Let us now turn to the more difficult case $\rho_k(\Phi_{f_\beta}) = 0.$
As before,
 define the operators $C_{\ell} = 1 - \psi(\ell) V_{\ell}$ and $O_{\ell} = 1 - \tau(\ell)\ell^{k+1} V_{\ell}$,
 \[ C_\ell, O_\ell\colon  \Symb_{\Gamma'}(\D_k(L))[f_\beta] \rightarrow \Symb_{\Gamma}(\D_k(L))[f_\beta]. \] A 
standard computation using (\ref{UlVl}) and (\ref{UlVlTl}) shows that for each factorization into positive integers $N/QR = st$, the
corresponding symbol
\[ \Phi_{f_\beta}^{s,t} :=  (\Phi_{f_\beta})_{| \prod_{\ell \mid s} O_\ell \prod_{\ell \mid t} C_\ell} \]
is a $U_\ell$-eigenvector for each $\ell \mid N/(QR)$, with eigenvalue $\tau(\ell)$ for $\ell \mid s$
and eigenvalue $\psi(\ell) \ell^{k+1}$ for $\ell \mid t$.
Since these systems of eigenvalues are distinct for distinct factorizations $N/QR = st$,
it follows that the vectors $\Phi_{f_\beta}^{s,t}$ are linearly independent as long as they are non-zero.
Furthermore, the vectors in $\{ \Phi_{f_\beta}^{s,t}: st = N/QR \}$
 are clearly linear combinations of the vectors in $\{(\Phi_{f_\beta})_{|V_{t}}: t \mid N/(QR) \}$ 
by definition. These two sets of have the same size, so the linear independence of the $\Phi_{f_\beta}^{s,t}$
implies that of the $(\Phi_{f_\beta})_{|V_{t}}$.

It therefore remains to prove that $\Phi_{f_\beta}^{s,t} \neq 0$.  We will prove the following statement:

\begin{quote}
{\em If $\Phi$ is a nonzero symbol in $\Symb_{\Gamma'}^{\Psi(-1)}(\D_k^\dagger(L))[f_\beta]$ such that $\rho_k(\Phi) = 0$,
 $\ell$ is a prime dividing $N/QR$, and $\lambda \in L$ is not of the form $\ell^\nu$ for an integer $\nu \ge k+2$, then
 $\Phi - \lambda \Phi_{|V_\ell} \neq 0$.}
\end{quote}

The proposition follows by inductively applying this claim for  all $\ell \mid N/(QR)$.
 To prove the claim, we first show that if $\Phi - \lambda \Phi_{|V_\ell} = 0$, then 
 $\Phi(\{\infty\}-\{c\})=0$ for any cusp $c$ in the $\Gamma$-equivalence class of $0$. Such a cusp may be written $c=\frac{a}{b}$ where $a,b$ are relatively prime integers, with $b$ coprime to $Np$ (and in particular $\ell \nmid b$). There exists  a positive integer $n$ such that $\ell^n \equiv 1 \pmod{b}$. We have $\Phi = \lambda^n \Phi_{|V_{\ell^n}}$.
Applying this equation to the divisor $\{\infty\}-\{a/b\}$, we get 
\begin{eqnarray*} 
 \Phi(\{\infty\}-\{a/b\})(f(z)) &=& \lambda^n \ell^{(-k-1)n} \Phi(\{\infty\} - \{\ell^n a/b\})( \ell^{nk} f(z/\ell^n)) \\ 
 & & \text{\ \ (using the def.\ of $V_\ell$ and (\ref{actpk})) } \\ 
 &=& (\lambda/\ell)^n \Phi(\{\infty\} - \{a/b\})(f(z/ \ell^n +m)), \\
& & \text{ where $m= \frac{a}{b}(\ell^n -1) \in \Z$, using $\Phi_{|\mat{1 & m \\ 0 & 1}}=\Phi$.}
\end{eqnarray*} 
  Therefore, if we set $\mu := \Phi(\{\infty\}-\{a/b\})$, then for all $f \in \D^\dag(L)$,
$$\mu(f(z)) = (\lambda/\ell)^n \mu(f(z/\ell^n+m)).$$
We claim that such a distribution $\mu$ is $0$. A change of variable $z \mapsto z + \ell^n m$ allows us to assume $m=0$. Hence
$\mu$ vanishes against all functions $g(z)$ of the form $f(z) - (\lambda/\ell)^n f(z/\ell^n)$. One 
sees easily that all functions $g \in \D^\dag(L)$ are of this form, unless $\lambda = \ell^\nu$, with $\nu$ an integer $\geq 1$, in which case we obtain all functions $g$ whose $(\nu-1)$-th derivative at $0$ is equal to $0$. On the other hand, since $\rho_k(\Phi)=0$,
$\mu$ also vanishes against any polynomial in $z$ of degree less than or equal to $k$. 

By our assumption on $\lambda$, if $\lambda = \ell^\nu$ with $\nu$ an integer $\geq 1$, then $\nu \le k+1$,
and hence any element $ \D^\dag(L)$ can be written as a sum of a polynomial of degree at most $k$ with a 
function whose $(\nu-1)$-th derivative at $0$ is equal to $0$.
It therefore follows that $\mu=0$.

We have thus proven that $\Phi(\{\infty\}-\{c\})=0$ for all 
$c \in \Gamma \cdot 0$.  We deduce that the restriction of $\Phi$ to 
$\Symb_{\Gamma,\Gamma \cdot 0}(\D_k^\dag(L))$ is $0$. By Lemma~\ref{kerrestriction}, $\Phi$ is thus a boundary modular symbol. Lemma~\ref{lemmabsymb0} below then shows that $\Phi=0$, which gives a contradiction.  This proves the italicized claim above, and hence the proposition.
\end{pf}

\begin{lemma} \label{lemmabsymb0} One has $\BSymb_{\Gamma}^{\psi(-1)}(\D_k^\dag(L))[f_\beta] = 0$.
\end{lemma} 
\begin{pf}
By Theorem~\ref{thetakisom}, it suffices 
to show that there is no Hecke eigenvector in
$\BSymb_{\Gamma}(\D_{-2-k}^\dag(L))$ with the same eigenvalues  as in Theorem~\ref{t:psidef} (in particular with $\iota$-eigenvalue $-\tau(-1)$).
By Theorem~\ref{basisboundover}, we know that a basis of  this eigenspace without the condition on the $\iota$-eigenvalue
 is given by the $\Phi_{-2-k,\tau,\psi,t}$ for $QRt|N$. 
 All these vectors have sign $\tau(-1)$, hence the lemma. 
\end{pf} 

\begin{prop} \label{dimfour} Let $S$ denote the set of divisors of $N/(QR)$.
The space of partial modular symbols
$\Symb_{\Gamma,C}^{\psi(-1)}(\D_k^\dag(L))[f_\beta]$
has dimension $|S|$ and is generated by the restriction of $(\Phi_{f_\beta})_{|V_{t}}$
 to that space for $t \in S$.
\end{prop}
\begin{pf} 
Let $x$ be the point of $\CC_{\Gamma,\BCM}$ corresponding to $f_\beta$.
By Theorem~\ref{t:smootheigen}, $x$ is a smooth point, and in a neighborhood of $x$ in $\CC_{\Gamma,\BCM}$, all classical points are old, but come from newforms of level $\Gamma' = \Gamma_1(QR) \cap \Gamma_0(p)$.

By Corollary~\ref{cordimxy}, the dimension of $\Symb_{\Gamma,C}^{\psi(-1)}(\D_k^\dag(L))[x]$ is the same as the dimension of $S_{k'+2}(\Gamma,L)[y]$ for $y$ a classical point close to $x$. By the above,
this is the dimension of the eigenspace in $S_{k'+2}(\Gamma,L)$ of the system of eigenvalues of a form of level $\Gamma'$.  Atkin-Lehner's theory of newforms gives that the dimension of this space is $|S|$. 
(To be precise, if $g(z)$ is a newform on $\Gamma_1(QR) \cap \Gamma_0(p)$, then the corresponding eigenspace for forms on $\Gamma$ has a basis $\{ g(tz): t \in S\}$.) This proves the first assertion.

The restriction map \[ \Symb_{\Gamma}^{\psi(-1)}(\D_k)[f_\beta] \rightarrow  \Symb_{\Gamma,C}^{\psi(-1)}(\D_k)[f_\beta] \]
is injective since by Lemma~\ref{kerrestriction} the kernel of this map is contained in 
 $\BSymb_{\Gamma}^{\psi(-1)}(\D_k)[f_\beta]$, which is $0$ by the Lemma~\ref{lemmabsymb0}.
 Hence the restrictions of the independent vectors (by Prop.~\ref{fourindependent})  $(\Phi_{f_\beta})_{|V_{t}}$ are still independent, hence a basis of $\Symb_{\Gamma,C}^{\psi(-1)}(\D_k^\dag(L))[f_\beta]$.
\end{pf}

To conclude the proof of Theorem~\ref{t:main} in the normal case, we write the symbol $\Phi_k^\crit$ of Corollary~\ref{corphikcrit}
in the  basis above: there exist constants $a_t \in L$ such that 
$$\Phi_k^\crit =  \sum_{t \in S} a_t (\Phi_{f_\beta})_{|V_t}.$$
Taking the Mellin transform, we get for all $\sigma \in \WW(\C_p)$ with $\sigma(-1)=\psi(-1)$:
$$L(\Phi_k^\crit,\sigma) = r(\sigma) L(\Phi_{f_\beta},\sigma),$$ where
\[ r(\sigma) = \sum_{t \in S} a_t t^{-1} \sigma(t)^{-1}. \]
From Corollary~\ref{corphikcrit} we obtain
\begin{equation} \label{e:algind}
\frac{L_p(\Phi_{f_\beta},\sigma)  }{  \log^{[k+1]}(\sigma) L_p(\psi,\sigma z) L_p(\tau,\sigma z^{-k}) \sigma^{-1}(R)}
= \frac{n(k, \sigma)}{r(\sigma)}.
\end{equation}
The left  side of (\ref{e:algind}) is independent of the primes dividing $N/(QR)$, 
and hence the right side must be unchanged if these primes are replaced by other primes.
However, the right side of (\ref{e:algind}) is a rational function in functions of the form $\sigma \mapsto \sigma(\ell)$.
Since any collection of functions of this form for distinct primes $\ell$ are an algebraically independent set of
functions on $\W^{\psi(-1)}$, the right side of
 (\ref{e:algind}) can be independent of the primes dividing $N/(QR)$ only if it is a constant.  From the defintion of $n(k, \sigma)$ and $r(\sigma)$, 
 this constant is clearly non-zero, and hence we obtain
 $${L_p(\Phi_{f_\beta},\sigma)} = {  \log^{[k+1]}(\sigma)  \sigma^{-1}(R) L_p(\psi,\sigma z) L_p(\tau,\sigma z^{-k})}$$
 up to a nonzero constant.

\section{Proof of the main theorem in the exceptional case}

In this section we consider the exceptional Eisenstein series $f=E_{2,\ell}$. 
Our proof of the main theorem for the normal case in Section~\ref{s:families} does not apply here because
Theorem~\ref{t:smootheigen} does not hold in the exceptional case (since the primes dividing $N/QR$ are necessarily bad),
and Prop.~\ref{dimfour} relies on Theorem~\ref{t:smootheigen}.  We therefore employ a more explicit method.

Let $p$ be a prime $\neq \ell$, and let $f_\beta = E_{2,\ell}(z)-E_{2,\ell}(pz)$ be the critical refinement of $f$. Let 
$\Phi_{f_\beta} \in \Symb^+_{\Gamma_0(p\ell)}(\D_0^\dag(\Q_p))$ the non-zero overconvergent modular symbol with the same  eigenvalues as $f_\beta$, i.e.\ such that
\begin{itemize}
\item[(i)] $T_{\ell'} \Phi_{f_\beta} = (1+\ell') \Phi_{f_\beta}$, for $\ell'$ prime, $\ell' \nmid \ell p$,
\item[(ii)] $U_p \Phi_{f_\beta} = p \Phi_{f_\beta}$.
\end{itemize}
The uniqueness of $\Phi_{f_\beta}$
up to multiplication by a non-zero scalar is guaranteed by the main result of \cite{Bcrit}.
We even have $\rho_0(\Phi_{f_\beta}) = \phi_{f_\beta}$ up to a non-zero scalar, where $\phi_{f_\beta} \in \Symb_{\Gamma(\ell p)}(\Q_p)$ is the classical boundary modular symbol attached to $f_\beta$, since the eigencurve is \'etale over the weight space at $f_\beta$ (cf. \cite{Bcourse}). 

In this section, we denote by $C_0$ the set of cusps $\Gamma_0(\ell) \cdot 0$ and $C_\infty$ the set of cusps $\Gamma_0(\ell) \cdot \infty$ (this supersedes the notation used in section~\ref{s:overconvergent}). Let us choose two integers $x$ and $y$ such that $\ell x - p y=1$ and let $w_{\ell}=\mat{\ell & y \\ \ell p & \ell x} \in \Gl_2(\Q)$. One has $\det \omega_\ell = \ell$. An easy computation shows that this matrix normalizes 
$\Gamma_0(\ell p)$, interchanges  $C_0$ and $C_\infty$, and satisfies $w_\ell^2 \in \ell \Sl_2(\Z)$.  For  an $S_0(p)$-module
$W$,
$w_\ell$ induces maps 
\begin{align*}
w_{\ell}\colon& \  \Symb_{\Gamma_0(\ell p)}(W) \rightarrow 
\Symb_{\Gamma_0(\ell p)}(W), \\ 
w_{\ell}\colon & \ \Symb_{\Gamma_0(\ell p),C_0}(W) \rightarrow 
\Symb_{\Gamma_0(\ell p),C_\infty}(W) 
\end{align*}
that are easily seen to be compatible with the Hecke operators $T_q$ for $q$ prime to $\ell p$, with $U_p$, and when $-\Id$ acts trivially on $W$, with the action of $\iota$. In particular, one sees that $w_\ell 
\Phi_{f_\beta}$ is a scalar times $\Phi_{f_\beta}$. To determine that scalar we use \cite[Lemma 5]{atkinlehner},
which states that $w_\ell(\phi_{f_\beta})+U_\ell(\phi_{f_\beta})$ has level $\Gamma_0(p)$, hence is $0$ since there
is no classical modular symbol (or form) on level $\Gamma_0(p)$ with the system of eigenvalues of $E_2^\crit$. Since 
$U_\ell(\phi_{f_\beta}) = \phi_{f_\beta}$, we get:
\begin{eqnarray} \label{e:wellphifb}  w_\ell \Phi_{f_\beta} = -\Phi_{f_\beta}. \end{eqnarray}

Let us choose an auxiliary prime $q$ not dividing $\ell p$. We define $\Gamma=\Gamma_0(q \ell p)$ and let 
$C_{q,0} = \Gamma_0(q \ell) \cdot 0$, $C_{q,\infty} = \Gamma_0(q\ell) \cdot \infty$ and 
$C_{q} =  C_{q,0} \coprod C_{q,\infty}$. These sets of cusps  fit into the following diagram, where 
arrows represent inclusions:
$$ 
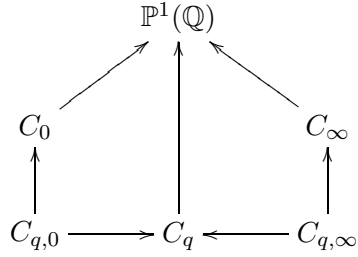
\xymatrix{  &  \PP^1(\Q) & \\ C_0 \ar[ur] & & C_\infty \ar[ul] \\  C_{q,0} \ar[u] \ar[r] & C_{q} \ar[uu]  & C_{q,\infty} \ar[u] \ar[l] }$$
For any sets of cusps $C \subset C'$ we denote by $\res_{C',C}$ the obvious restriction map on partial modular symbols, and just $\res_C$ when $C'$ is the full set of cusps $\PP^1(\Q)$. 
We also have operators 
\begin{align*}
V_q \colon & \  \Symb_{\Gamma_0(\ell p),C_0} \rightarrow  \Symb_{\Gamma_0(\ell q p),C_0}, \\
V_q\colon & \ \Symb_{\Gamma_0(\ell p),C_\infty} \rightarrow  \Symb_{\Gamma_0(\ell q p),C_\infty}.
\end{align*}

Using the construction of Darmon--Dasgupta \cite{DD}, one proves the following result:
\begin{prop} \label{p:phiDD} There exist non-zero partial overconvergent modular symbols 
\begin{align*}
\Phi_{\infty} \in& \  \Symb^+_{\Gamma_0(\ell p),C_\infty}(\D_0^\dag(\Q_p)) \\
\Phi_{0} \in& \ \Symb^+_{\Gamma_0(\ell p),C_0}(\D_0^\dag(\Q_p)) \\
\Phi_{q} \in & \ \Symb^+_{\Gamma,C_{q}}(\D_0^\dag(\Q_p))
\end{align*}
such that 
\begin{itemize}
\item[(i)] $T_{\ell'} \Phi_{\infty} = (1+\ell') \Phi_{\infty}$ and $T_{\ell'} \Phi_{0} = (1+\ell') \Phi_{0}$ for $\ell' \nmid \ell p$.
\item[(ii)] $U_p \Phi_{\infty} = p \Phi_{\infty}$ and $U_p \Phi_{0} = p \Phi_{0}.$
\item[(iii)]  $w_{\ell} \Phi_{0} = - \Phi_{\infty}$.
\item[(iv)] One has $\res_{C_{q},C_{q,0}} \Phi_{q} =  \res_{C_0,C_{q,0}} (1- q V_{q}) \Phi_0$ and similarly with $0$ replaced by 
$\infty$.
\item[(v)] The $p$-adic Mellin transform of the distribution $\Phi_{q}(\{\infty\}-\{0\})$ is 
$$L_p(\Phi_{q}) (\sigma) =  \log_p^{[1]}(\sigma) (1-\sigma^{-1}(q))(1-\sigma^{-1}(\ell)) \zeta_p(\sigma z) \zeta_p(\sigma) \quad \text{ if } \sigma(-1)=1.$$
\end{itemize}
\end{prop}
\begin{pf} 
Consider the usual classical modular form \[ E_{k+2}^\ord(z) = E_{k+2}(z)-p^{k+1} E_{k+2}(z) \] of level $\Gamma_0(p)$ for $k \geq 0$ an even integer,  
 with $E_{k+2}(z)$ the usual Eisenstein series of level $1$ and weight $k+2$ (which is not a modular form for $k+2=2$, even though $E_{2}^\ord$ is).
To this form we add three others
\begin{eqnarray*}
 F_{k+2}^\infty &=& (1-V_\ell) E_{k+2}^\ord \in M_{k+2}(\Gamma_0(\ell p)) \\ 
F_{k+2}^0 &=& (1- \ell^{k+2} V_\ell) E_{k+2}^\ord  \in M_{k+2}(\Gamma_0(\ell p))  \\
F_{k+2}^{q} &=& (1-q^{k+2} V_{q})(1 - V_\ell) E_{k+2}^\ord  \in M_{k+2}(\Gamma) 
\end{eqnarray*}
As is well known, $E_{k+2}^\ord$ is part of a family indexed by the weight space, and thus in particular makes sense for any $k \in \Z$, as do $F_{k+2}^\infty, F_{k+2}^0$ and $F_{k+2}^{q}$.
The crucial point for the proof is the following ``numerical coincidence" that occurs at $k=-2$: 
$$F_{0}^\infty = F_{0}^0.$$
Therefore,
\begin{eqnarray} \label{numcoin} F_0^{q} = (1-V_{q}) F_0^0 = (1-V_{q}) F_0^\infty.\end{eqnarray}

It is easy to see (see the comments above Definition~\ref{d:Ndef}) that $F_{k+2}^\infty$ is $C_\infty$-cuspidal, $F_{k+2}^0$ is $C_0$-cuspidal, and $F_{k+2}^{q}$ is $C_{q}$-cuspidal.  
We can therefore attach classical partial modular symbols to those three modular forms (cf. Definition~\ref{d:defpartial}), and lift them  (since they are ordinary at $p$) to partial overconvergent modular symbols 
\begin{align*}
\Phi_{k+2}^\infty \in & \ \Symb_{\Gamma_0(\ell p),C_\infty}(\D^\dag_k(\Q_p)), \\
\Phi_{k+2}^0 
\in & \ \Symb_{\Gamma_0(\ell p),C_0}(\D^\dag_k(\Q_p)), \\
 \Phi_{k+2}^{q} \in & \
\Symb_{\Gamma,C_{q}}(\D^\dag_k(\Q_p)).
\end{align*}
As $k$ varies, each of these three symbols is   part of a family of ordinary overconvergent modular symbols (cf. \cite{DD}) over the weight space and thus makes sense for  any $k$ in $\Z$. In particular for $k=-2$ we get 
\begin{align*}
\Phi_{0}^\infty \in & \  \Symb_{\Gamma_0(\ell p),C_\infty}(\D^\dag_{-2}(\Q_p)), \\
 \Phi_{0}^0 \in & \ \Symb_{\Gamma_0(\ell p),C_0}(\D^\dag_{-2}(\Q_p)), \\
 \Phi^{q}_0 \in  & \ \Symb_{\Gamma,C_{q}}(\D^\dag_{-2}(\Q_p)).
 \end{align*}
  Finally we define $\Phi_\infty$, $\Phi_0$ and $\Phi_{q}$ appearing in the statement of the proposition as the image by $\Theta_0$ (cf.~(\ref{fundexsymb})) of these modular symbols.

 It is  easy to compute by interpolation the Hecke eigenvalues of $\Phi_{0}^\infty$ and $\Phi_{0}^0$, and thus 
 those of $\Phi_\infty$ and $\Phi_0$ which proves (i) and (ii). 
 
 A direct computation using the 
 convergent series defining $E_{k+2}$ shows that
 $w_{\ell} F_{k+2}^0 = - \ell^{-1} F_{k+2}^\infty$ for $k \geq 0$.  By interpolation, a similar relation holds for the attached modular symbols, including for $k=-2$: $w_\ell \Phi_{0}^0 = - \ell^{-1} \Phi_{0}^\infty$. Applying $\Theta_0$ and taking into account 
 $w_\ell \Theta_0 = \ell \Theta_0 w_\ell$,    we get (iii). 
 
 The point (iv) is a translation of the numerical coincidence (\ref{numcoin}): one gets from (\ref{numcoin}) that
 $$\res_{C_{q},C_{q,0}} \Phi^{q}_0 =  \res_{C_0,C_{q,0}} (1-  V_{q}) \Phi^0_0$$ and similarly with $0$ replaced by $\infty$.
  After applying $\Theta_0$ we get (iv), noting that $\Theta_0 V_{q} = q V_{q} \Theta_0$.

  The point (v) follows by computing the $p$-adic $L$-function attached to the form $F_{k+2}^{q}$, as in Prop.~\ref{p:lfk},
 interpolating this to $k=-2$ as in Theorem~\ref{t:psidef} and applying $\Theta_0$ as in Corollary~\ref{corphikcrit}.
\end{pf}

\begin{prop} \label{p:phires0i} Up to multiplying $\Phi_{f_\beta}$ by a non-zero scalar, we have
\begin{eqnarray} \label{e:resco} \res_{C_0} \Phi_{f_\beta} &=& \Phi_0 \\
\label{e:resci} \res_{C_\infty} \Phi_{f_\beta} &=& \Phi_\infty
\end{eqnarray}
N.B.\ the  multiplication by a single non-zero scalar is enough to imply (\ref{e:resco}) and (\ref{e:resci}) simultaneously.
 \end{prop}
\begin{pf}
Since the eigencurve of full modular symbols and tame level $\ell$ is smooth (even étale over the weight space) at the point $E_{2,\ell}^\crit=f_\beta$, as is the eigencurve of partial $C_0$-modular symbols, the eigenspaces 
\[ \Symb^+_{\Gamma_0(p\ell)}(\D_0^\dag(\Q_p))[f_\beta] \text{ and  } \Symb^+_{\Gamma_0(\ell p),C_0}(\D_0^\dag(\Q_p))[f_\beta]\] have dimension $1$.
Therefore $\res_{C_0}  \Phi_{f_\beta} = \alpha_0 \Phi_0$ for some $\alpha_0 \in \Q_p$. If $\alpha_0=0$, then $\Phi_{f_\beta}$ is a boundary modular symbol by Lemma~\ref{kerrestriction}, which contradicts Remark~\ref{r:signboundary}
since $\Phi_{f_\beta}$ has sign $+1$. So $\alpha_0 \neq 0$.

Similarly, one proves that $\res_{C_\infty} \Phi_{f_\beta} = \alpha_\infty \Phi_{\infty}$ with $\alpha_\infty \neq 0$.
Finally, we recall that $w_\ell \Phi_{f_\beta} = - \Phi_{f_\beta}$ while
by (iv) of  Proposition~\ref{p:phiDD}, $w_\ell \Phi_0 = - \Phi_\infty$.
The compatibility between $w_\ell$ and the restriction maps $\res_{C_0}$ and $\res_{C_\infty}$
 implies that $\alpha_0 = \alpha_\infty$.
\end{pf}

By  Proposition~\ref{p:phires0i} and Proposition~\ref{p:phiDD}(iv) one has
\begin{align} \res_{C_{q},C_{q,0}} (\res_{C_{q}} (1-q V_{q}) \Phi_{f_\beta}) =& \ \res_{C_{q},C_{q,0}} \Phi_{q} \label{e:res1}
\\
 \res_{C_{q},C_{q,\infty}} (\res_{C_{q}} (1- q V_{q}) \Phi_{f_\beta}) =& \ \res_{C_{q},C_{q,\infty}}  \Phi_{q} \label{e:res2}
\end{align}
We are thus in situation to apply the following general lemma:
\begin{lemma} Let $W$ be a $\Gamma$-module such that $W^\Gamma=0$. Then the map
$$\res_{C_{q},C_{q,0}} \times \res_{C_{q},C_{q,\infty}} :
\Symb_{\Gamma,C_{q}}(W) \rightarrow \Symb_{\Gamma,C_{q,0}}(W) \times  \Symb_{\Gamma,C_{q,\infty}}(W)$$ is injective.
\end{lemma}
\begin{pf} Let $\phi$ be in the kernel of the given map. We need to show that $\phi(\{a\}-\{b\})=0$ when $a \in C_{q,0}$ and $b \in C_{q,\infty}$. Since $\phi$ is in the kernel, $\phi(\{a\} -  \{0\})=0$ and $\phi(\{\infty\}-\{b\})=0$, and $\phi(\{0\}-\{\infty\}) \in W^\Gamma = 0$ by hypothesis. The result follows.
\end{pf}
The hypothesis $(\D_0^\dag)^\Gamma = 0$  of the lemma  
 is satisfied by \cite[Prop. 3.1]{stevenspollack2}.
 Therefore by (\ref{e:res1}) and (\ref{e:res2}), we obtain:
 $$ \res_{C_{q}} (1-qV_{q}) \Phi_{f_\beta} = \Phi_{q}.$$
Evaluating these modular symbols at the divisor $\{0\}-\{\infty\}$, which is in $C_{q}$, we get using
Prop.~\ref{p:phiDD}(v):
$$ (1-\sigma^{-1}(q)) L_p(f_\beta,\sigma) =  \log_p^{[1]}(\sigma)(1- \sigma^{-1}(q)) (1-\sigma^{-1}(\ell)) \zeta_p(\sigma z) \zeta_p(\sigma) $$
 for $\sigma$ such that $\sigma(-1)=1$.
Hence, canceling  the factor $1- \sigma^{-1}(q)$, we obtain Theorem~\ref{t:main} in the exceptional case.

\par \bigskip
\end{document}